\documentclass[12pt]{article}

\usepackage{lmodern}
\usepackage{geometry}
\usepackage{fullpage}
\usepackage{microtype}
\usepackage[utf8]{inputenc}
\usepackage[T1]{fontenc}
\usepackage{setspace}

\usepackage{enumerate}

\usepackage[round]{natbib}

\usepackage{xcolor}
\usepackage{hyperref}
\usepackage[title]{appendix}

\usepackage{textcomp}

\usepackage{graphicx}
\usepackage{floatrow}
\usepackage{caption, subcaption}
\usepackage{booktabs,threeparttable,array}

\usepackage{amsmath, amssymb, amsthm}
\theoremstyle{plain}
\newtheorem{theorem}{Theorem}
\newtheorem{lemma}{Lemma}
\newtheorem{proposition}[theorem]{Proposition}
\newtheorem{corollary}{Corollary}
\theoremstyle{definition}

\newtheorem{assumption}{Assumption}

\theoremstyle{remark}

\newenvironment{assumptionprime}[1]
  {\renewcommand{\theassumption}{\ref*{#1}$'$}%
   \addtocounter{assumption}{-1}%
   \begin{assumption}}
  {\end{assumption}}

\newenvironment{assumptionpprime}[1]
  {\renewcommand{\theassumption}{\ref*{#1}$''$}%
   \addtocounter{assumption}{-1}%
   \begin{assumption}}
  {\end{assumption}}

\newcommand{\E}{\mathsf{E}}
\newcommand{\var}{\mathrm{var}}
\newcommand{\cov}{\mathrm{cov}}
\newcommand{\cor}{\mathrm{cor}}
\newcommand{\pr}{\mathrm{P}}
\newcommand{\dd}{\mathrm{d}}
\newcommand{\veps}{\varepsilon}
\newcommand{\tr}{\operatorname{tr}}
\newcommand{\diag}{\operatorname{diag}}
\newcommand{\tsp}{\mathsf{T}}
\newcommand{\rN}{\mathrm{N}}
\newcommand{\R}[1]{\mathbb{R}^{#1}}
\DeclareMathOperator*{\argmax}{arg\,max}
\newcommand{\weak}{\rightsquigarrow}
\newcommand{\Proj}{\Pi}

\newcommand{\zipcovnsim}{100{,}000}
\newcommand{\zipvalidmaxdev}{0.013}
\newcommand{\zipnaivemin}{0.792}
\newcommand{\zipnaiveonemin}{0.763}
\newcommand{\zipwaldmin}{0.947}
\newcommand{\zipwaldmax}{0.996}
\newcommand{\ziplrtmin}{0.949}
\newcommand{\ziplrtmax}{0.976}
\newcommand{\ziplfmin}{0.940}
\newcommand{\zippreddiff}{0.013}
\newcommand{\zippredonediff}{0.018}
\newcommand{\zippimin}{0.946}
\newcommand{\zippimax}{0.963}
\newcommand{\zipboxpmin}{0.913}
\newcommand{\zipboxbmin}{0.955}
\newcommand{\zipboxbmax}{0.966}
\newcommand{\zipboxmmin}{0.948}
\newcommand{\zipboxmmax}{0.956}
\newcommand{\zipmbarmin}{2.12}
\newcommand{\zipmbarmax}{2.19}

\begin{document}

\title{Near-boundary asymptotics for profile score statistics}

\author{Matias Shedden \quad Karl Oskar Ekvall \\
{\normalsize Department of Statistics, University of Florida} \\
{\tt \normalsize sheddenmatias@ufl.edu \quad k.ekvall@ufl.edu }}
\date{October 24, 2026}

\maketitle

\onehalfspacing

\begin{abstract}
	We establish asymptotic results for profile score statistics under
sequences of parameters that may tend to the boundary of the parameter set.
When the likelihood is profiled over a set with boundary points, the
distribution of the statistic can depend, even asymptotically, on a distance to
the boundary that cannot be estimated consistently. We address this in three
ways. First, we give conditions under which the nuisance score has a root on an
extended parameter set, and the profile score statistic formed with that
root converges to a chi-squared distribution. Second, we show that the
efficient score statistic converges to a chi-squared distribution under a
separate condition that can hold even for estimators constrained to the original
set. Both lead to confidence regions with asymptotically correct coverage
probability, uniformly on sets that may contain boundary points, and the
conditions allow the scores for the interest and nuisance parameters to grow at
different rates or become asymptotically collinear. Third, we quantify the
coverage error of applying the profile score statistic without an extended
set. We approximate the distribution of the statistic and bound the error by
a quantity that decays like a normal tail in the standardized distance to the
boundary. The approximation also gives one-sided tests with least
favorable critical values whose asymptotic size is at most the nominal level. We
illustrate the results in a zero-inflated Poisson model and a one-way random
effects model.

\end{abstract}

\vspace{1em}
\vspace{1em}
\newpage
\section{Introduction}  \label{sec:introduction}

Inference near the boundary of the parameter set is needed in many
settings, for instance when boundary points correspond to simpler models
or a necessarily nonnegative parameter is small. Examples include
near-negligible
zero-inflation in a model for counts, random effect variances near zero or
correlations near one in a mixed model, a rare component in a mixture model, and
weak overdispersion in a negative binomial model. Sometimes the parameter of
interest is near the boundary, and sometimes a nuisance parameter is. In either
case, standard asymptotic approximations are inaccurate. For
every sample size, there are interior parameter values at which coverage
probabilities of conventional confidence regions are substantially different from nominal.

Let $\theta = (\psi, \lambda) \in \Theta$ be a parameter, where $\psi \in \R{d_1}$ is the
parameter of interest and $\lambda \in \R{d_2}$ the nuisance parameter. Let
$\ell_n(\theta) = \ell_n(\psi, \lambda)$ be a log-likelihood,
$S_{n\psi}(\theta)$ and $S_{n\lambda}(\theta)$ the components of the score
$S_n(\theta)$, $\mathcal{I}_n(\theta)$ the expected information matrix with blocks
$\mathcal{I}_{n\psi\psi}$, $\mathcal{I}_{n\psi\lambda} =
\mathcal{I}_{n\lambda\psi}^\tsp$ and $\mathcal{I}_{n\lambda\lambda}$, and $\bar{\mathcal{I}}_n =
\mathcal{I}_{n\psi\psi} -
\mathcal{I}_{n\psi\lambda}\mathcal{I}_{n\lambda\lambda}^{-1}\mathcal{I}_{n\lambda\psi}$
the efficient information. The subscript $n$ indexes a sequence of models and
corresponding quantities, and it need not be a sample size. For an estimator
$\hat\lambda_n(\psi)$ of the nuisance parameter at a fixed $\psi$, let
$\hat\theta_n = \{\psi, \hat\lambda_n(\psi)\}$. We consider profile score
statistics
\begin{equation} \label{eq:Rn}
  \hat R_n(\psi) = \bar{\mathcal{I}}_n(\hat\theta_n)^{-1/2} S_{n\psi}(\hat\theta_n),
\end{equation}
for different choices of the nuisance parameter estimator. The usual profile score is obtained by taking a partial maximizer, $\hat{\lambda}_n(\psi) \in \argmax_{\lambda \in
\Lambda(\psi)} \ell_n(\psi, \lambda)$, where $\Lambda(\psi) = \{\lambda : (\psi,
\lambda) \in \Theta\}$. We will see that this choice often leads to unreliable inference when $\lambda$ is near the boundary of $\Lambda(\psi)$, which motivates other choices.

To be more specific, consider the confidence region
\begin{equation*}
C_n(\alpha) = \{\psi : \|\hat R_n(\psi)\|^2 \leq q_{d_1, 1-\alpha}\},
\end{equation*}
where $q_{d_1, 1-\alpha}$ is the $1-\alpha$ quantile of a $\chi^2_{d_1}$
distribution. Our main results give conditions that ensure asymptotically
correct coverage of $C_n(\alpha)$, and of regions from other score-based
statistics and nuisance estimators discussed below, uniformly on sets $K
\subseteq \Theta$ that may contain boundary points. That is,
\begin{equation} \label{eq:coverage}
  \lim_{n \to \infty} \sup_{\theta \in K} \bigl| \pr_\theta\{\psi \in C_n(\alpha)\} - (1 - \alpha) \bigr| = 0.
\end{equation}
The main complications are that $\psi$ and $\lambda$ can each be near or on the
boundary of their respective parameter sets. The latter is the more challenging
because, while the statistic is evaluated at the null value of $\psi$, $\lambda$
is estimated. Natural choices of $\hat\lambda_n(\psi)$, including the partial
maximizer above, can be on the boundary with nonvanishing probability, in which
case $\hat R_n(\psi)$ is not asymptotically $\chi_{d_1}^2$-distributed.

Score statistics have previously proved useful for inference near the
boundary without nuisance parameters. The asymptotic distributions of likelihood
ratio and Wald statistics are nonstandard at the boundary, with or without
nuisance parameters \citep{chernoff1954distribution, self1987asymptotic,
geyer1994asymptotics, andrews1999estimation, andrews2001testing,
kopylev2011asymptotic, bertinellisalucci2026asymptotic, ekvall2026asymptotics},
and they depend on the limit of a standardized distance from a sequence of true
parameters to the boundary \citep{feder1968distribution, andrews2010asymptotic,
cox2022generalized}. By contrast, the score standardized by the expected
information often has the same limit along every sequence, which gives
\eqref{eq:coverage} without nuisance parameters \citep{bottai2003confidence,
ekvall2022confidence}. However, those results do not apply or extend easily to
settings with nuisance parameters. In fact, $C_n(\alpha)$ with the usual profile
score often has worse coverage than the Wald and likelihood ratio intervals it
is meant to improve on. For example, in simulations from a zero-inflated Poisson
model with the zero-inflation probability at or near zero, $C_n(0.05)$ based on
the usual profile score has coverage as low as \zipnaivemin{}, against at least
\zipwaldmin{} for the Wald and likelihood ratio intervals
(Section~\ref{sec:zip}).

\citet{andrews2009hybrid, andrews2010asymptotic} allow nuisance parameters near
the boundary and ensure correct asymptotic size with a critical value that is
least favorable over the limit of a standardized distance from the true
parameter to the boundary, which cannot be estimated consistently. The test
therefore rejects with probability below the nominal level at most parameter
values. \citet{mccloskey2017bonferronibased} instead takes the least favorable
value only over a confidence set for the standardized distance.
\citet{ketz2018subvector} assumes an extremum estimator over $\Theta$ is
$n^{1/2}$-consistent along sequences and bases a test on an unconstrained
one-step update that is asymptotically normal near and at the boundary. The
critical values are conditional on an asymptotically sufficient statistic for
$\lambda$ rather than least favorable. \citet{ketz2025short} note that inverting
the tests of \citet{andrews2009hybrid}, \citet{mccloskey2017bonferronibased} and
\citet{ketz2018subvector} to obtain confidence regions is often computationally
prohibitive, especially when $\lambda$ has more than one or two dimensions.
\citet{cox2026simple} gives a uniformly valid confidence interval when $d_1 = 1$
and $\Theta$ is defined by a single inequality. The method requires an estimator
of $\theta$ that is $n^{1/2}$-consistent and asymptotically normal.

Whether \eqref{eq:coverage} holds depends on the distribution of $\hat
R_n(\psi_n)$ along every sequence of true parameters $\theta_n = (\psi_n,
\lambda_n)$ in $K$, including ones approaching the boundary. When
$\hat\lambda_n(\psi)$ is a partial maximizer over $\Lambda(\psi)$, that
distribution depends on a standardized distance from $\lambda_n$ to the boundary
of $\Lambda(\psi_n)$, specifically the distance in the metric $(\lambda,
\lambda') \mapsto \|\mathcal{I}_{n\lambda\lambda}(\theta_n)^{1/2}(\lambda -
\lambda')\|$. Along sequences where the distance is bounded,
$\hat\lambda_n(\psi_n)$ can be on the boundary with probability bounded away
from zero, $\hat R_n(\psi_n)$ need not be asymptotically normal, and the
distance cannot be estimated consistently. We address this in three ways. None of the three requires independent observations or that different
blocks or eigenvalues of $\mathcal{I}_n(\theta_n)$ grow at the same rate. In particular, when $n$ is a sample size, $n^{-1}\mathcal{I}_n(\theta_n)$ need not converge. The
true $\theta_n$ may approach the boundary at any rate, so there are no restrictions on the standardized distance to the boundary. Except in the third
approach, correlations between elements of $S_{n\psi}(\theta_n)$ and $S_{n\lambda}(\theta_n)$ may tend to one in absolute value.

First, we let $\hat\lambda_n(\psi)$ be a root of the nuisance score on an
extended parameter set, which may depend on $n$. The statistic then behaves as
if the true parameter were an interior point, provided the standardized distance
from $\lambda_n$ to the boundary of the extended set diverges. This distance in
general differs from the corresponding distance in the original set. While
extended parameter sets have been considered previously
\citep{self1987asymptotic, feng1992statistical}, formal conditions ensuring
asymptotic normality of $\hat R_n(\psi_n)$ along sequences have not been
established. This is addressed by Theorems \ref{thm:prof} and \ref{thm:roots}.
Proposition \ref{prop:box} extends Theorem \ref{thm:prof} to simultaneous
intervals for several coordinates of $\theta$.

Second, for settings where the first approach is unavailable, we
consider the efficient score statistic $\hat R^{\mathrm{e}}_n$
\citep{neyman1959optimal}, which replaces $S_{n\psi}(\hat\theta_n)$ in
\eqref{eq:Rn} by $S_{n\psi}(\hat\theta_n) - 
\mathcal{I}_{n\psi\lambda}(\hat\theta_n)\mathcal{I}_{n\lambda\lambda}(\hat\theta_n)^{-1}
S_{n\lambda}(\hat\theta_n)$.
Theorem \ref{thm:eff} shows the statistic is asymptotically standard normal
under conditions that are different from those needed for $\hat{R}_n$, and in
particular are compatible with estimators that can be on the boundary with
probability bounded away from zero. \citet{fan2023wald} use a statistic that
coincides with $\hat R^{\mathrm{e}}_n$ in some cases and obtain tests with
nominal asymptotic size
when the constraints are linear, all parameters are estimated at
rate $n^{1/2}$, and the limiting information is nonsingular. Any consistent
estimator of the information, observed or expected, can then be used.
Theorem \ref{thm:eff} requires none of these conditions, and its statistic uses
the expected information. Standardizing by the observed information instead can
give a limit that depends on the sequence \citep{bottai2003confidence,
ekvall2022confidence}. Theorems \ref{thm:prof} and \ref{thm:eff} give
\eqref{eq:coverage} with chi-squared critical values, because each statistic has
the same limit along every sequence. This makes implementation straightforward.

Third, we provide approximations to the distribution of $\hat R_n(\psi_n)$
when $\hat{\lambda}_n(\psi)$ is a partial maximizer over $\Lambda(\psi)$, in
which case \eqref{eq:coverage} does not hold in general, as noted above. In
fact, in this setting $\hat{R}_n(\psi_n)$ need not converge in distribution at
all. Theorem \ref{thm:approx} shows that $\hat R_n(\psi_n)$ and 
a variable $R^\ast_n$, built from a standard normal vector and the true $\theta_n$,
give asymptotically the same probability to every convex
set, uniformly over the convex sets and over $K \subseteq \Theta$ that may contain boundary points, for arbitrary but fixed $d_1$
and $d_2$. When $d_1 = d_2 = 1$ and $\lambda_n$ is near the lower endpoint of
$\Lambda(\psi_n)$ only, the approximating variable is
\[
  R^\ast_n = Z_\psi + \kappa_n(1 - \kappa_n^2)^{-1/2}\min(Z_\lambda + h_n, 0),
\]
where $Z_\psi$ and $Z_\lambda$ are independent standard normal, $\kappa_n$ is
the correlation between $S_{n\psi}(\theta_n)$ and $S_{n\lambda}(\theta_n)$, and
$h_n$ is the standardized distance from $\lambda_n$ to the lower endpoint of
$\Lambda(\psi_n)$. The correlation can be consistently estimated but $h_n$
cannot. The coverage error of $C_n(\alpha)$ decays like a normal tail in $h_n$
(Proposition \ref{prop:bound}). Taking $\kappa_n$ at its estimate and $h_n$ at
its least favorable value, as in the plug-in size correction of
\citet[Section~3.3]{andrews2009hybrid}, gives one-sided tests with asymptotic
size at most the nominal level, uniformly over $K$ and in the level (Proposition
\ref{prop:lf}). Near the boundary, the one-sided test against alternatives of
sign opposite to $\kappa_n$ is more powerful than the corresponding test based
on $\hat R^{\mathrm{e}}_n$, and far from it less powerful. The other one-sided
test is never more powerful.

Section \ref{sec:examples} illustrates the results in a zero-inflated
Poisson model, with $\psi$ the Poisson mean and $\lambda$ the zero-inflation
probability, and in a one-way random effects model, with $\psi$ the error
variance and $\lambda$ the between-unit variance. In both, the distance from
$\lambda_n$ to the boundary of the extended set diverges. With the partial maximizer over
that set, \eqref{eq:coverage} holds on every compact subset of $\Theta$ in the
first model and on all of $\Theta$ in the second. A
companion paper \citep{shedden2026score}, extending \citet{ekvall2026uniform}
and \citet{zhang2025fast}, verifies the conditions of Theorem \ref{thm:prof} in
linear mixed models with several variance components.


\section{Uniform asymptotics for profile and efficient scores}
\label{sec:general_theory}

\subsection{Model and statistics} \label{sec:model}

The parameter set $\Theta$ is a subset of $\R{d}$, with $d = d_1 + d_2$.
For each $\lambda \in \R{d_2}$, define $\Psi(\lambda) = \{\psi \in \R{d_1} :
(\psi, \lambda) \in \Theta\}$, the analogue of $\Lambda(\psi)$ of Section
\ref{sec:introduction}. For
$n \in \{1, 2, \dots\}$, $\ell_n:\Theta \to \R{}$ is a log-likelihood function
based on a family of distributions indexed by $\theta$, and $\theta_n = (\psi_n,
\lambda_n)$ is the true parameter. We write $\E_\theta$ for expectations under
the distribution indexed by $\theta$, and similarly $\var_\theta$,
$\cov_\theta$, and $\pr_\theta$; we omit the subscript when it is $\theta_n$,
unless it needs emphasis. The score of Section \ref{sec:introduction} is
$S_n(\theta) = \partial\ell_n(\theta)/\partial\theta$, and the information
matrix is $\mathcal{I}_n(\theta) = \cov_\theta\{S_n(\theta)\}$. The efficient
information of Section \ref{sec:introduction} is $\bar{\mathcal{I}}_n(\theta) =
\mathcal{I}_{n\psi\psi}(\theta) - \beta_n(\theta)\mathcal{I}_{n\lambda\psi}(\theta)$,
where $\beta_n(\theta) =
\mathcal{I}_{n\psi\lambda}(\theta)\mathcal{I}_{n\lambda\lambda}(\theta)^{-1}$ is
the coefficient in the population regression of $S_{n\psi}(\theta)$ on
$S_{n\lambda}(\theta)$. Let $\mathcal{J}_n(\theta) =
-\partial S_n(\theta)/\partial\theta^\tsp$ be the observed information, with the
same block notation.

A key object in the following is the efficient score for $\psi$,
\begin{equation} \label{eq:eff_score}
  U_n(\theta) = S_{n\psi}(\theta) - \beta_n(\theta_n) S_{n\lambda}(\theta),
\end{equation}
in which the coefficient is evaluated at $\theta_n$ for every argument $\theta$.
The efficient information satisfies $\bar{\mathcal{I}}_n(\theta_n) = \cov\{U_n(\theta_n)\}$.
A tilde on a score denotes the score
standardized by the square root of its covariance matrix at $\theta_n$,
whatever the argument of the score: $\tilde S_n(\theta) =
\mathcal{I}_n(\theta_n)^{-1/2} S_n(\theta)$, $\tilde S_{n\lambda}(\theta) =
\mathcal{I}_{n\lambda\lambda}(\theta_n)^{-1/2} S_{n\lambda}(\theta)$, and
$\tilde U_n(\theta) = \bar{\mathcal{I}}_n(\theta_n)^{-1/2} U_n(\theta)$. In
general, $\tilde S_{n \lambda}$ is not a block of $\tilde S_{n}$ unless
$\mathcal{I}_n(\theta_n)$ is block-diagonal.

The nuisance estimators $\hat{\lambda}_n(\psi)$ of Section
\ref{sec:introduction} are statistics only when a value of $\psi$ is
hypothesized, but for simplicity we call them estimators. For
example, $\hat{\lambda}_n(\psi)$ may be a root of the nuisance score over an
extended parameter set, a partial maximizer of the log-likelihood over
$\Lambda(\psi)$,
or a plug-in estimator. In general, $\hat\lambda_n(\psi)$ need not be in
$\Lambda(\psi)$, provided the functions taking it as an argument can be extended
appropriately; Section \ref{sec:roots} makes this precise.

An efficient score statistic for $\psi$ takes $\hat{\theta}_n = \{\psi, \hat\lambda_n(\psi)\}$ in place of $\theta_n$ and $\theta$ in \eqref{eq:eff_score}. Specifically,
\begin{equation} \label{eq:eff_stat}
   \hat{R}^{\mathrm{e}}_n(\psi) = \bar{\mathcal{I}}_n(\hat{\theta}_n)^{-1/2}
   \{S_{n\psi}(\hat{\theta}_n) - \beta_n(\hat{\theta}_n) S_{n\lambda}(\hat{\theta}_n)\}.
\end{equation}
The statistic $\|\hat{R}^{\mathrm{e}}_n(\psi)\|^2$ specializes to the test
statistic of \citet{neyman1959optimal} when evaluated at a $n^{1/2}$-consistent
nuisance estimator. The profile score statistic \eqref{eq:Rn} drops the
last term of \eqref{eq:eff_stat}. When $d_1
= 1$, $\hat{R}_n$ is the signed root of the score statistic $\|\hat{R}_n\|^2$. The statistics differ in general but coincide if $\hat\lambda_n(\psi)$ is a
root of the nuisance score, $S_{n\lambda}(\hat{\theta}_n) = 0$. They also
coincide if $\psi$ and $\lambda$ are orthogonal in the Fisher information sense,
that is, if $\mathcal{I}_{n\psi\lambda}(\theta) = 0$ for every $\theta$.

The region $C^{\mathrm{e}}_n(\alpha)$ is $C_n(\alpha)$ with
$\hat{R}^{\mathrm{e}}_n$ in place of $\hat{R}_n$. The following lemma reduces \eqref{eq:coverage} to conditions along sequences of
true parameters; its proof is in Appendix \ref{sec:chi2_proofs}.

\begin{lemma} \label{lem:unif}
Fix $K \subseteq \Theta$ and $c \in \R{}$, and let $g_n : \Theta \to \R{}$
for $n \in \{1, 2, \dots\}$. It holds that
\begin{enumerate}
  \item[(i)~~] if $g_n(\theta_n) \rightarrow c$ for every sequence $(\theta_n)$
    in $K$, then $\sup_{\theta \in K}|g_n(\theta) - c| \rightarrow 0$; and
  \item[(ii)~] if $\limsup_{n \rightarrow \infty} g_n(\theta_n) \leq c$ for every
    sequence $(\theta_n)$ in $K$, then $\limsup_{n \rightarrow \infty}
    \sup_{\theta \in K} g_n(\theta) \leq c$.
\end{enumerate}
\end{lemma}

We say that a condition holds on a set $K \subseteq \Theta$ if it holds
along every sequence $(\theta_n)$ in $K$. A consequence of this lemma is that
if $\hat{R}_n(\psi_n) \rightsquigarrow \rN_{d_1}(0, I_{d_1})$ on $K$, then
\eqref{eq:coverage} holds. Indeed, $\|\hat{R}_n(\psi_n)\|^2 \rightsquigarrow
\chi^2_{d_1}$ by the continuous mapping theorem, and because the $\chi^2_{d_1}$ distribution function is continuous,
the convergence in distribution implies $\pr_{\theta_n}\{\psi_n \in C_n(\alpha)\}
\rightarrow 1 - \alpha$ for every $\alpha \in (0, 1)$, so Lemma
\ref{lem:unif}(i) applies with $g_n(\theta) = \pr_\theta\{\psi \in
C_n(\alpha)\}$ and $c = 1 - \alpha$.

\subsection{Chi-squared limits and uniform coverage} \label{sec:limits}

We first state high-level conditions under which the asymptotic properties of
the test statistics hold. These are later verified in models and given primitive
conditions that may be easier to verify in applications. The conditions are
stated for an arbitrary sequence $(\theta_n)$ of true parameters.

\begin{assumption} \label{assm:information}
 The information matrix $\mathcal{I}_n(\theta_n)$ is invertible for every $n$.
\end{assumption}

Assumption \ref{assm:information} ensures that the efficient score and the
efficient information at $\theta_n$ are well defined, and that the latter is
positive definite. Indeed, $\mathcal{I}_n(\theta_n)$ is an invertible
covariance matrix, hence positive definite, and a symmetric matrix is positive
definite if and only if a principal block and its Schur complement are
\citep[Theorem~7.7.7]{horn2012matrix}. We use this fact repeatedly.
Notably, Assumption \ref{assm:information} does not require the
eigenvalues and diagonal elements of $\mathcal{I}_n(\theta_n)$ to grow at the
same rate. For example, the correlation between the scores for interest and
nuisance parameters may tend to one in absolute value, making the efficient information
arbitrarily small relative to $\mathcal{I}_{n\psi\psi}(\theta_n)$.

\begin{assumption} \label{assm:clt}
  The standardized score is asymptotically normal in the sense that $\tilde S_n(\theta_n) \rightsquigarrow \rN_d(0, I_d)$.
\end{assumption}

We do not assume $\E\{S_n(\theta_n)\} = 0$, but Assumption \ref{assm:clt}
implies $\E\{\tilde S_n(\theta_n)\} \to 0$: the standardized score is tight and,
centered, has covariance $I_d$, so its mean and second moments are bounded, and
uniform integrability \citep[Theorem~2.20]{vandervaart1998asymptotic} gives the
limit. Without nuisance parameters, Assumption \ref{assm:clt} on a set $K
\subseteq \Theta$ suffices, by Lemma \ref{lem:unif}, for score-based confidence
regions to have asymptotically correct coverage probability uniformly on $K$.
The information matrix need not be invertible at every point of interest:
if $K$ is dense in the set of interest, a continuity argument can extend the
coverage to points where the information matrix is singular
\citep{ekvall2022confidence}.

At $\psi = \psi_n$, write $\hat\lambda_n = \hat\lambda_n(\psi_n)$, so that
$\hat\theta_n = (\psi_n, \hat\lambda_n)$. By \eqref{eq:Rn}, \eqref{eq:eff_score} and \eqref{eq:eff_stat},
\begin{equation} \label{eq:numerators}
\begin{aligned}
  \hat R_n(\psi_n)
  &= \bar{\mathcal{I}}_n(\hat\theta_n)^{-1/2}U_n(\hat\theta_n) + \bar{\mathcal{I}}_n(\hat\theta_n)^{-1/2}\beta_n(\theta_n)S_{n\lambda}(\hat\theta_n), \\
 \hat R^{\mathrm{e}}_n(\psi_n)
  &= \bar{\mathcal{I}}_n(\hat\theta_n)^{-1/2}U_n(\hat\theta_n) - \bar{\mathcal{I}}_n(\hat\theta_n)^{-1/2}\{\beta_n(\hat\theta_n) - \beta_n(\theta_n)\}
  S_{n\lambda}(\hat\theta_n).
\end{aligned}
\end{equation}
The first term, common to both statistics, is the efficient score at
$\hat\theta_n$ standardized by the efficient information at $\hat\theta_n$. With
$\theta_n$ in place of $\hat\theta_n$, this term would be the standardized
efficient score $\tilde U_n(\theta_n)$. The following assumptions ensure that
this replacement is asymptotically negligible, in the efficient information and
in the efficient score, respectively. Whether they hold depends on both the
structure of these functions and the choice of $\hat\lambda_n$. The trailing
terms in \eqref{eq:numerators}, which differ between the two statistics, are
dealt with in Assumptions \ref{assm:nuis_est} and
\ref{assm:coef_stable} below. For a matrix, $\|\cdot\|$ is the spectral norm.

\begin{assumption} \label{assm:eff_info_stable}
  The efficient information is stable under plugging in $\hat\lambda_n$ in the sense that
  \[
    \| \bar{\mathcal{I}}_n(\theta_n)^{1/2} \bar{\mathcal{I}}_n(\hat{\theta}_n)^{-1} \bar{\mathcal{I}}_n(\theta_n)^{1/2} - I_{d_1} \| = o_p(1).
  \]
\end{assumption}

The inverses in Assumption \ref{assm:eff_info_stable} need not exist
for every outcome, so the converging quantity need not be a random variable.
For such quantities, convergence means that the quantity exists with
probability tending to one and that the statement holds for every extension
to a random variable. In particular, it is part of Assumption
\ref{assm:eff_info_stable} that the inverses exist with probability tending
to one.

\begin{assumption} \label{assm:linear_approx}
  The efficient score is stable under plugging in $\hat\lambda_n$ in the sense that
  \[
    \tilde U_n(\hat\theta_n) - \tilde U_n(\theta_n) = o_p(1).
  \]
\end{assumption}

If $\ell_n$ is twice differentiable, then differentiating \eqref{eq:eff_score} gives
\[
  \partial U_n(\psi_n, \lambda)/\partial\lambda^\tsp
  = \beta_n(\theta_n)\mathcal{J}_{n\lambda\lambda}(\psi_n, \lambda) - \mathcal{J}_{n\psi\lambda}(\psi_n, \lambda),
\]
which vanishes at $\lambda = \lambda_n$ when $\mathcal{J}_n$ is replaced by
$\mathcal{I}_n$, by the definition of $\beta_n$. To first order, the left-hand
side of Assumption \ref{assm:linear_approx} is therefore the deviation of the
observed information from $\mathcal{I}_n(\theta_n)$ multiplied by the nuisance
estimation error. Consequently, Assumption \ref{assm:linear_approx} holds, for
example, when the observed information is close to $\mathcal{I}_n(\theta_n)$
near $\theta_n$ and $\hat\lambda_n - \lambda_n$ is of the order of
$\mathcal{I}_{n\lambda\lambda}(\theta_n)^{-1/2}$, as in Theorem \ref{thm:roots}.

The following lemma says that, under these assumptions, the first, shared term
in \eqref{eq:numerators} is asymptotically standard normal.

\begin{lemma} \label{lem:chi2}
If Assumptions \ref{assm:information}--\ref{assm:linear_approx} hold along a
sequence $(\theta_n)$ of true parameters, then
$\bar{\mathcal{I}}_n(\hat\theta_n)^{-1/2}U_n(\hat\theta_n) \rightsquigarrow
\rN_{d_1}(0, I_{d_1})$.
\end{lemma}

\begin{proof}
We first show $\tilde U_n(\theta_n) \rightsquigarrow \rN_{d_1}(0, I_{d_1})$;
it suffices to show every subsequence has a further one along which
it holds \citep[Theorem 2.6]{billingsley1999convergence}. Let
\[
  \Gamma_n = \bar{\mathcal{I}}_n(\theta_n)^{-1/2}[I_{d_1},
  -\beta_n(\theta_n)]\mathcal{I}_n(\theta_n)^{1/2}.
\]
Because $U_n(\theta_n) = [I_{d_1}, -\beta_n(\theta_n)]S_n(\theta_n)$, we have
$\tilde U_n(\theta_n) = \Gamma_n \tilde S_n(\theta_n)$, and because
$\cov\{\tilde S_n(\theta_n)\} = I_d$, $\Gamma_n\Gamma_n^\tsp = \cov\{\tilde
U_n(\theta_n)\} = I_{d_1}$ by the definition of the efficient information.
Thus, $(\Gamma_n)$ is contained in a compact set, and every subsequence has a
further subsequence along which $\Gamma_n \to \Gamma$ with $\Gamma\Gamma^\tsp =
I_{d_1}$. Along it, $\{\tilde S_n(\theta_n), \Gamma_n\} \rightsquigarrow (Z,
\Gamma)$ with $Z \sim \rN_d(0, I_d)$ by Assumption \ref{assm:clt}, and hence
$\tilde U_n(\theta_n) \rightsquigarrow \Gamma Z \sim \rN_{d_1}(0, I_{d_1})$ by
the continuous mapping theorem.

Now let $M_n =
\bar{\mathcal{I}}_n(\hat\theta_n)^{-1/2}\bar{\mathcal{I}}_n(\theta_n)^{1/2}$,
so that
\[
  \bar{\mathcal{I}}_n(\hat\theta_n)^{-1/2}U_n(\hat\theta_n) =
  M_n[\tilde U_n(\theta_n) + \{\tilde U_n(\hat\theta_n) - \tilde U_n(\theta_n)\}].
\]
Assumption \ref{assm:eff_info_stable} says $M_n^\tsp M_n - I_{d_1} = o_p(1)$,
so by Lemma \ref{lem:sqrt}, with probability
tending to one $\bar{\mathcal{I}}_n(\hat\theta_n)$ is positive definite and
$\|M_n - I_{d_1}\|^2 \leq 3 d_1 \|M_n^\tsp M_n - I_{d_1}\|$; hence
$M_n - I_{d_1} = o_p(1)$, and the term in braces is $o_p(1)$ by Assumption
\ref{assm:linear_approx}. Slutsky's theorem gives the conclusion.
\end{proof}

Let
\begin{equation}  \label{eq:An}
  \tilde\beta_n(\theta) = \bar{\mathcal{I}}_n(\theta)^{-1/2} \beta_n(\theta)
  \mathcal{I}_{n\lambda\lambda}(\theta)^{1/2} \in \R{d_1 \times d_2}.
\end{equation}
The matrix $\tilde\beta_n(\theta)$ is the coefficient in the (population)
regression of $\bar{\mathcal{I}}_n(\theta)^{-1/2}S_{n\psi}(\theta)$ on
$\mathcal{I}_{n\lambda\lambda}(\theta)^{-1/2}S_{n\lambda}(\theta)$ under
$\theta$. By \eqref{eq:eff_score}, $\bar{\mathcal{I}}_n(\theta_n)^{-1/2}
S_{n\psi}(\theta) = \tilde U_n(\theta) + \tilde\beta_n(\theta_n) \tilde
S_{n\lambda}(\theta)$ for every $\theta$. Hence $\tilde\beta_n(\theta_n)\tilde
S_{n\lambda}(\hat\theta_n)$ is the second term in the first line of
\eqref{eq:numerators} with $\theta_n$ in place of $\hat\theta_n$ in the
efficient information, and the following assumption makes it negligible.

\begin{assumption} \label{assm:nuis_est}
  The $\hat\lambda_n$ used in $\hat{R}_n$ satisfies
  \[
    \tilde\beta_n(\theta_n) \tilde S_{n\lambda}(\hat\theta_n) = o_p(1).
  \]
\end{assumption}
If $\mathcal{I}_{n\psi\lambda}(\theta_n) = 0$, then $\tilde\beta_n(\theta_n) = 0$ and Assumption
\ref{assm:nuis_est} holds for every $\hat\lambda_n$. Otherwise, it holds if
$\hat\lambda_n$ is a root of the nuisance score, exactly or with probability
tending to one. A partial maximizer over $\Lambda(\psi_n)$ need not be a root
when $\lambda_n$ is near the boundary, but a root over an extended nuisance set
can often be used instead (Theorem~\ref{thm:roots}). In Section \ref{sec:irreg}
we treat settings where the nuisance estimator can be on the boundary and
Assumption \ref{assm:nuis_est} fails, in which case the profile score statistic need
not be asymptotically chi-squared.

\begin{theorem} \label{thm:prof}
If Assumptions \ref{assm:information}--\ref{assm:nuis_est} hold on a set
$K \subseteq \Theta$, then \eqref{eq:coverage} holds.
\end{theorem}

\begin{proof}
Fix a sequence $(\theta_n)$ in $K$ and let $M_n =
\bar{\mathcal{I}}_n(\hat\theta_n)^{-1/2}\bar{\mathcal{I}}_n(\theta_n)^{1/2}$.
By the first identity in \eqref{eq:numerators} and by \eqref{eq:An},
\[
  \hat R_n(\psi_n) = \bar{\mathcal{I}}_n(\hat\theta_n)^{-1/2}U_n(\hat\theta_n)
  + M_n\tilde\beta_n(\theta_n)\tilde S_{n\lambda}(\hat\theta_n).
\]
The first term is asymptotically standard normal by Lemma \ref{lem:chi2}. The
second is $o_p(1)$ because $\tilde\beta_n(\theta_n)\tilde
S_{n\lambda}(\hat\theta_n) = o_p(1)$ by Assumption \ref{assm:nuis_est} and
$\|M_n\|^2 = \|M_n^\tsp M_n\| = 1 + o_p(1)$ by Assumption
\ref{assm:eff_info_stable}. Hence
$\hat R_n(\psi_n) \rightsquigarrow \rN_{d_1}(0, I_{d_1})$. The sequence was
arbitrary, so \eqref{eq:coverage} holds, as shown after Lemma \ref{lem:unif}.
This completes the proof.
\end{proof}

The following assumption can hold in settings where Assumption
\ref{assm:nuis_est} fails because the nuisance estimator can be on the boundary
of $\Lambda(\psi_n)$.

\begin{assumptionprime}{assm:nuis_est} \label{assm:coef_stable}
  The $\hat\lambda_n$ used in $\hat{R}^{\mathrm{e}}_n$ satisfies
  \[
    \bar{\mathcal{I}}_n(\theta_n)^{-1/2}\{\beta_n(\hat\theta_n) -
    \beta_n(\theta_n)\} S_{n\lambda}(\hat\theta_n) = o_p(1).
  \]
\end{assumptionprime}

Neither of Assumptions \ref{assm:nuis_est} and \ref{assm:coef_stable} implies
the other. Often,
$\mathcal{I}_{n\psi\lambda}(\theta_n) \neq 0$, so that $\beta_n(\theta_n) \neq
0$, and $\beta_n(\hat\theta_n) - \beta_n(\theta_n)$ is of smaller
order than $\beta_n(\theta_n)$. Assumption \ref{assm:coef_stable} then allows a
larger nuisance score. However, unlike Assumption \ref{assm:nuis_est}, Assumption \ref{assm:coef_stable} is not automatically satisfied when $\mathcal{I}_{n\psi\lambda}(\theta_n) = 0$.

\begin{theorem} \label{thm:eff}
If Assumptions \ref{assm:information}--\ref{assm:linear_approx} and
\ref{assm:coef_stable} hold on a set $K \subseteq \Theta$, then
\eqref{eq:coverage} holds with $C_n^{\mathrm{e}}(\alpha)$
in place of $C_n(\alpha)$.
\end{theorem}

\begin{proof}
Fix a sequence $(\theta_n)$ in $K$. The proof of Theorem \ref{thm:prof}
applies with the second identity in \eqref{eq:numerators} in place of the
first, which replaces $\tilde\beta_n(\theta_n)\tilde S_{n\lambda}(\hat\theta_n)$
by $-\bar{\mathcal{I}}_n(\theta_n)^{-1/2}\{\beta_n(\hat\theta_n) -
\beta_n(\theta_n)\}S_{n\lambda}(\hat\theta_n)$, and with Assumption
\ref{assm:coef_stable} in place of Assumption \ref{assm:nuis_est}. Hence
$\hat{R}^{\mathrm{e}}_n(\psi_n) \rightsquigarrow \rN_{d_1}(0, I_{d_1})$. The
sequence was arbitrary, so the argument after Lemma \ref{lem:unif}, with
$\hat{R}^{\mathrm{e}}_n$ and $C^{\mathrm{e}}_n(\alpha)$ in place of $\hat{R}_n$
and $C_n(\alpha)$, gives the conclusion. This completes the proof.
\end{proof}

Theorem \ref{thm:eff} can be used with partial maximizers over $\Lambda(\psi)$, which
need no extended set. Every quantity in $\hat{R}^{\mathrm{e}}_n$ is then evaluated
in $\Theta$. This is useful when the score or information has no convenient
continuation outside $\Theta$ or when software returns only estimates in
$\Theta$, as is common for generalized linear mixed models and variance components models, for example.

\subsection{Roots on an extended set} \label{sec:roots}

We now give conditions on the observed information that ensure the
nuisance score has a root on an extended set, so that Assumption
\ref{assm:nuis_est} holds, and under which Assumptions
\ref{assm:eff_info_stable} and \ref{assm:linear_approx} hold for every
estimator satisfying
\begin{equation} \label{eq:info_bounded}
  \|\mathcal{I}_{n\lambda\lambda}(\theta_n)^{1/2}(\hat\lambda_n - \lambda_n)\| = O_p(1).
\end{equation}

For each $\psi$ and $n$, let $\tilde\Lambda_n(\psi) \supseteq
\Lambda(\psi)$ be a set to which $S_n(\psi, \cdot)$ and $\mathcal{I}_n(\psi,
\cdot)$ extend, the latter with symmetric values; $\ell_n$ need not be defined,
and $(\psi, \lambda)$ need not index a distribution, at the added points. In the
zero-inflated Poisson model of Section \ref{sec:zip}, the extended set is a
deterministic interval of zero-inflation probabilities that contains negative
values and is the same for each $n$. In the linear mixed models of
\citet{shedden2026score}, it is the set of covariance parameters for which the
covariance matrix of the data is positive definite, whether or not that of the
random effects is. That extended set is different for each $n$. For $M > 0$, let
\[
  B_n^M = \{\lambda \in \tilde\Lambda_n(\psi_n) : \|\mathcal{I}_{n\lambda\lambda}(\theta_n)^{1/2}(\lambda -
  \lambda_n)\| \leq M\},
\]
and
\begin{equation} \label{eq:rdist}
  r_n = \inf_{\lambda \notin \Lambda(\psi_n)}
  \|\mathcal{I}_{n\lambda\lambda}(\theta_n)^{1/2}(\lambda - \lambda_n)\|, \qquad
  \tilde r_n = \inf_{\lambda \notin \tilde\Lambda_n(\psi_n)}
  \|\mathcal{I}_{n\lambda\lambda}(\theta_n)^{1/2}(\lambda - \lambda_n)\|,
\end{equation}
which are the distances from $\lambda_n$ to the complements of the original and
extended nuisance sets, in the metric of $\mathcal{I}_{n\lambda\lambda}(\theta_n)$.

The hypothesis $\tilde r_n \to \infty$ of Theorem \ref{thm:roots} below says the extended set is large on the
information scale. In the examples of Section \ref{sec:examples}, $\tilde r_n$
diverges. In both, the true parameter can be on the boundary of the original set,
whereas the boundary of the extended set is outside the model. Section \ref{sec:irreg} treats the regime in which $r_n$ stays
bounded.

Write $\mathcal{J}_{n\cdot\lambda}(\theta) = -\partial S_n(\theta)/\partial\lambda^\tsp$
for the last $d_2$ columns of $\mathcal{J}_n(\theta)$, which exist where $S_n(\psi,
\cdot)$ is differentiable, and $\mathcal{I}_{n\cdot\lambda}$ for the same
columns of $\mathcal{I}_n$. The dot distinguishes these columns from the
$\lambda$ rows of $\mathcal{J}_n$. The rows contain $\mathcal{J}_{n\lambda\psi}(\theta) = -\partial
S_{n\lambda}(\theta)/\partial\psi^\tsp$, which need not exist under the
hypotheses of Theorem \ref{thm:roots} below, since those concern $S_n(\psi_n,
\cdot)$ only. Consider the conditions
\begin{equation} \label{eq:info_unif}
\begin{aligned}
  \sup_{\lambda \in B_n^M} \bigl\| \mathcal{I}_n(\theta_n)^{-1/2}\{\mathcal{J}_{n\cdot\lambda}(\psi_n, \lambda) -
  \mathcal{I}_{n\cdot\lambda}(\theta_n)\}\mathcal{I}_{n\lambda\lambda}(\theta_n)^{-1/2} \bigr\| &= o_p(1), \\
  \sup_{\lambda \in B_n^M} \bigl\| \mathcal{I}_n(\theta_n)^{-1/2}\{\mathcal{I}_n(\psi_n, \lambda) -
  \mathcal{I}_n(\theta_n)\}\mathcal{I}_n(\theta_n)^{-1/2} \bigr\| &= o(1).
\end{aligned}
\end{equation}
The first line says the observed
information is close to $\mathcal{I}_n(\theta_n)$, and the second that
$\mathcal{I}_n(\psi_n, \cdot)$ is stable, for $\lambda$ within a fixed distance
of $\lambda_n$ in that metric.

\begin{theorem} \label{thm:roots}
Suppose that, along a sequence $(\theta_n)$ in $\Theta$, Assumptions
\ref{assm:information} and \ref{assm:clt} hold, $S_n(\psi_n, \cdot)$ is
continuously differentiable on the interior of $\tilde\Lambda_n(\psi_n)$,
$\tilde r_n \to \infty$, and \eqref{eq:info_unif} holds for every $M > 0$.
Then
\begin{enumerate}
  \item[(i)~~] Assumptions \ref{assm:eff_info_stable} and
  \ref{assm:linear_approx} hold for every $\hat\lambda_n \in
  \tilde\Lambda_n(\psi_n)$ satisfying \eqref{eq:info_bounded};
  \item[(ii)~] there is a $\hat\lambda_n \in \tilde\Lambda_n(\psi_n)$
  satisfying \eqref{eq:info_bounded} that, with probability tending to one,
  is a root of $S_{n\lambda}(\psi_n, \cdot)$ with
  $\mathcal{J}_{n\lambda\lambda}(\psi_n, \hat\lambda_n)$ positive definite, and any
  root satisfying \eqref{eq:info_bounded} equals $\hat\lambda_n$ with
  probability tending to one.
\end{enumerate}
\end{theorem}

The proof is in Appendix \ref{sec:chi2_proofs}. A root of the nuisance
score satisfies Assumption \ref{assm:nuis_est} trivially, so by part (i) and
Theorem \ref{thm:prof}, $C_n(\alpha)$ formed with a root satisfying
\eqref{eq:info_bounded} has asymptotically correct coverage uniformly on any
$K$ on which the hypotheses hold, which may contain boundary points of
$\Theta$, with the estimator computed as if there were no boundary. Part (ii) says the estimator is well defined: the root exists and is the only one within a bounded distance of $\lambda_n$ on the information scale. If $\ell_n(\psi_n, \cdot)$ extends to
$\tilde\Lambda_n(\psi_n)$ with gradient $S_{n\lambda}(\psi_n, \cdot)$ it is a
strict local maximizer, so a partial maximizer over $\tilde\Lambda_n(\psi_n)$ satisfying
\eqref{eq:info_bounded} is that root. Part (i) also applies to estimators
that are not roots, such as a partial maximizer over $\Lambda(\psi_n)
\subseteq \tilde\Lambda_n(\psi_n)$. Section \ref{sec:zip}
uses it for partial maximizers over both nuisance sets.

\subsection{Simultaneous intervals} \label{sec:simul}

In practice, it is often desired to construct simultaneous intervals for
several coordinates of $\theta$ with simultaneous coverage, that is, a box. Bonferroni adjustment ensures validity but is conservative when the correlations in the inverse
information are large, and more so with more coordinates. For $j \in \{1, \dots,
d\}$, let $\hat{R}_{n,j}(\theta_j)$ be $\hat{R}_n$ of \eqref{eq:Rn} with the
$j$th coordinate $\theta_j$ as the parameter of interest and the other
coordinates as nuisance parameters, estimated by
$\hat\theta_{n,-j}(\theta_j)$. Let $\Theta_j$ be the set of $j$th coordinates
of points in $\Theta$. Let $\Omega_n(\theta)$ be the correlation matrix of
$\mathcal{I}_n(\theta)^{-1}$, let $\hat{\Omega}_n$ be an estimator of
$\Omega_n(\theta_n)$, and for a correlation matrix $\Omega$ let $m_\alpha(\Omega)$ be the
$(1 - \alpha)$-quantile of $\|\Omega^{1/2}Z\|_\infty$, $Z \sim \rN_d(0, I_d)$.
When $d = 2$, the off-diagonal element of $\Omega_n(\theta)$ is minus the
correlation between the two scores at $\theta$. The box is
\[
  C_n^{\mathrm{s}}(\alpha) = \prod_{j = 1}^d \bigl\{\theta_j \in \Theta_j : |\hat{R}_{n,j}(\theta_j)|
  \leq m_\alpha(\hat{\Omega}_n)\bigr\}.
\]

\begin{proposition} \label{prop:box}
If, on a set $K \subseteq \Theta$, Assumptions
\ref{assm:information} and \ref{assm:clt} hold, Assumptions
\ref{assm:eff_info_stable}--\ref{assm:nuis_est} hold with, for every $j$, the
$j$th coordinate as the parameter of interest and $\hat\theta_{n,-j}(\theta_{n,j})$
as the estimator, and $\hat{\Omega}_n - \Omega_n(\theta_n) = o_p(1)$, then
\eqref{eq:coverage} holds with $C_n^{\mathrm{s}}(\alpha)$ and $\theta$ in place of
$C_n(\alpha)$ and $\psi$.
\end{proposition}

The proof is in Appendix \ref{sec:chi2_proofs}. With $z_p$ the
$p$-quantile of the standard normal distribution, $z_{1-\alpha/2} \leq
m_\alpha(\Omega) \leq z_{1-\alpha/(2d)}$ for every $\Omega$, and the upper
bound is the Bonferroni critical value.

\section[Approximating distributions and least favorable critical values]{Approximating distributions and least favorable critical values}
\label{sec:irreg}

We now consider profile score confidence regions formed with a
partial maximizer $\hat\lambda_n(\psi) \in \argmax_{\lambda \in \Lambda(\psi)}
\ell_n(\psi, \lambda)$. Such estimators are common in practice, but they often
do not satisfy Assumption \ref{assm:nuis_est}, so Theorem \ref{thm:prof} does
not apply. More precisely, when $(r_n)$ of \eqref{eq:rdist} is bounded,
$\hat\lambda_n$ can be on the boundary of $\Lambda(\psi_n)$ with probability
not tending to zero, the nuisance score there need not be zero, and hence
Assumption \ref{assm:nuis_est} need not hold.

The main result of this section, Theorem \ref{thm:approx}, approximates the
distribution of the profile score statistic $\hat{R}_n$ uniformly over sets of
true parameters, whether or not $(r_n)$ is bounded. The approximating
distribution is asymptotically normal when $r_n \to \infty$ and can be far from
normal otherwise. It can be computed by one-dimensional integration when
$\psi$ and $\lambda$ are scalar, and comparing it with the normal approximation that
holds for the efficient score statistic of Theorem \ref{thm:eff} and for the
extended set roots of Section \ref{sec:roots} shows where a
partial maximizer over $\Lambda(\psi)$ gains power and where it loses.

Define the local nuisance set
\[
  \mathcal{H}_n = \mathcal{I}_{n\lambda\lambda}(\theta_n)^{1/2}\bigl\{\Lambda(\psi_n) - \lambda_n\bigr\}
  \subseteq \R{d_2},
\]
and let $\Proj_n$ be the Euclidean projection onto it; the
projection is well defined under Assumptions \ref{assm:information} and
\ref{assm:proj}.

\begin{assumptionpprime}{assm:nuis_est}\label{assm:proj} The set
$\Lambda(\psi_n)$ is closed and convex for every $n$, and $\hat\lambda_n$
satisfies, with $\tilde\beta_n$ as in \eqref{eq:An},
\[
  \tilde\beta_n(\theta_n) \tilde S_{n\lambda}(\hat\theta_n) = \tilde\beta_n(\theta_n)\bigl[ \tilde S_{n\lambda}(\theta_n) - \Proj_n\{\tilde S_{n\lambda}(\theta_n)\}\bigr] + o_p(1).
\]
\end{assumptionpprime}

Assumption \ref{assm:proj} is similar to Assumption \ref{assm:nuis_est} but with
an additional term on the right-hand side. The additional term is needed because
$\hat{\lambda}_n$ can be on the boundary with probability not tending to zero.
If $r_n \to \infty$, then Assumption \ref{assm:proj} reduces to Assumption
\ref{assm:nuis_est}. Indeed, $\mathcal{H}_n$ contains the ball of radius $r_n$
about the origin, and $\tilde S_{n\lambda}(\theta_n)$ is tight by Assumption
\ref{assm:clt}, so if $r_n \to \infty$ then $\tilde S_{n\lambda}(\theta_n) \in
\mathcal{H}_n$, and hence $\Proj_n\{\tilde S_{n\lambda}(\theta_n)\} = \tilde
S_{n\lambda}(\theta_n)$, with probability tending to one.

The following lemma gives sufficient conditions for Assumption
\ref{assm:proj} that are sometimes easier to verify. It uses Assumption
\ref{assm:bddcorr} below, and its proof is in Appendix \ref{sec:irreg_proofs}.

\begin{lemma}\label{lem:primitive}
If $\Lambda(\psi_n)$ is closed and convex for every $n$,
Assumption~\ref{assm:bddcorr} holds, and
\begin{enumerate}
  \item[(i)~~] $\displaystyle
    \mathcal{I}_{n\lambda\lambda}(\theta_n)^{1/2}(\hat\lambda_n - \lambda_n)
    = \Proj_n\{\tilde S_{n\lambda}(\theta_n)\} + o_p(1)$,
  \item[(ii)~] $\displaystyle
    \tilde S_{n\lambda}(\hat\theta_n) - \tilde S_{n\lambda}(\theta_n)
    = -\mathcal{I}_{n\lambda\lambda}(\theta_n)^{1/2}(\hat\lambda_n - \lambda_n) + o_p(1)$,
\end{enumerate}
then Assumption~\ref{assm:proj} holds.
\end{lemma}

Condition (i) says $\hat\lambda_n$ behaves as the maximum likelihood estimator
in a normal shift experiment with parameter set $\mathcal{H}_n$. It is a
near-boundary analogue of the result that a constrained maximum likelihood
estimator is asymptotically the projection of the score onto a cone
\citep{geyer1994asymptotics}. Here, $\mathcal{H}_n$ and
$\tilde\beta_n(\theta_n)$ need not converge, so there is no single limiting
experiment to pass to. Condition (ii) says the nuisance score expands to first
order about $\theta_n$ with $\mathcal{I}_{n\lambda\lambda}(\theta_n)$ in place
of the observed information.

Since $\hat{R}_n(\psi_n)$ need not converge in distribution, we compare its
distribution with that of an approximating variable at each $n$, and show
that the two give asymptotically the same probability to every convex set,
uniformly over sets of true parameters.
The approximating variable is
\begin{equation}
  R^\ast_n = Z_\psi + \tilde\beta_n(\theta_n)\{Z_\lambda - \Proj_n(Z_\lambda)\} \in \R{d_1},
  \label{eq:Tstar}
\end{equation}
where $Z = (Z_\psi, Z_\lambda) \sim \rN_d(0, I_d)$. The second term is zero when
$Z_\lambda \in \mathcal{H}_n$; intuitively, it is otherwise the deviation from normality
due to the estimator being on the boundary of $\Lambda(\psi_n)$. The variable $R^\ast_n$
approximates $\hat R_n(\psi_n)$ because, by Assumptions
\ref{assm:linear_approx} and \ref{assm:proj},
\begin{equation}
  \bar{\mathcal{I}}_n(\theta_n)^{-1/2} S_{n\psi}(\hat\theta_n)
  = \tilde U_n(\theta_n) + \tilde\beta_n(\theta_n)[\tilde S_{n\lambda}(\theta_n) - \Proj_n\{\tilde S_{n\lambda}(\theta_n)\}] + o_p(1).
  \label{eq:Vnlin}
\end{equation}
The left-hand side is asymptotically equivalent to $\hat R_n(\psi_n)$ by Assumption \ref{assm:eff_info_stable} and Lemma \ref{lem:sqrt}. The
variable $R^\ast_n$ is the right-hand side with $Z_\psi$ and $Z_\lambda$ in
place of the standardized scores $\tilde U_n(\theta_n)$ and $\tilde
S_{n\lambda}(\theta_n)$, which are uncorrelated with identity covariance
matrices and, by Assumption \ref{assm:clt}, jointly asymptotically standard
normal. Thus, $R_n^\ast$ depends on $n$ only through the nonstochastic
quantities $\mathcal{H}_n$ and $\tilde\beta_n(\theta_n)$, that is, through the
geometry of the nuisance set and the information matrix.

The results in this section require an asymptotic bound on the linear dependence
between the scores.
When $d_1 = d_2 = 1$, let $\kappa_n(\theta) = \cor_\theta\{S_{n\psi}(\theta),
S_{n\lambda}(\theta)\}$; otherwise let it be the
largest canonical correlation between $S_{n\psi}(\theta)$ and $S_{n\lambda}(\theta)$
\citep[Section 12.2]{anderson2003introduction}. Write $\kappa_n = \kappa_n(\theta_n)$.
By Lemma \ref{lem:Akappa}, $\|\tilde\beta_n(\theta_n)\|^2 =
\kappa_n^2/(1 - \kappa_n^2)$ and, when $d_1 = d_2 = 1$, $\tilde\beta_n(\theta_n) = \kappa_n/(1 -
\kappa_n^2)^{1/2}$.

\begin{assumption}\label{assm:bddcorr}
The matrices $\tilde\beta_n(\theta_n)$ of \eqref{eq:An} satisfy $\limsup_{n \to \infty}
\|\tilde\beta_n(\theta_n)\| < \infty$ or, equivalently, $\limsup_{n \to \infty} |\kappa_n| < 1$.
\end{assumption}

To state the first result of the section, let $\mathcal{C}$ be the family of convex Borel subsets of $\R{d_1}$.
For $\theta \in \Theta$, let
$R^\ast_n(\theta)$ be $R^\ast_n$ with $\theta$ as the true parameter, so that
$R^\ast_n = R^\ast_n(\theta_n)$.

\begin{theorem}\label{thm:approx}
If Assumptions~\ref{assm:information}--\ref{assm:linear_approx},
\ref{assm:proj} and \ref{assm:bddcorr} hold on a set $K \subseteq \Theta$,
then
\[
  \sup_{\theta \in K}\ \sup_{C \in \mathcal{C}}\ \bigl| \pr_\theta\{\hat{R}_n(\psi) \in C\}
  - \pr\{R^\ast_n(\theta) \in C\} \bigr| \to 0 .
\]
\end{theorem}
The balls $\{x \in \R{d_1} : \|x\|^2 \leq q_{d_1, 1-\alpha}\}$, $\alpha \in (0, 1)$, are in $\mathcal{C}$,
so Theorem \ref{thm:approx} has the following corollary.

\begin{corollary} \label{cor:naive}
If the assumptions of Theorem \ref{thm:approx} hold, then
\[
  \sup_{\theta \in K}\ \sup_{\alpha \in (0, 1)}\ \bigl| \pr_{\theta}\{\psi \in C_n(\alpha)\}
  - \pr\{\|R^\ast_n(\theta)\|^2 \leq q_{d_1, 1-\alpha}\} \bigr| \to 0 .
\]
\end{corollary}

Since $R_n^*$ is essentially $\hat{R}_n$ with the score replaced by its normal limit, the corollary suggests any substantial coverage distortion of $C_n(\alpha)$ is due to the parameter set and information matrix geometry, not due to non-normality of the score. Thus, analyzing the behavior of $R_n^*$ provides insight into the coverage properties of $C_n(\alpha)$. The following result gives an upper bound on how different from normal $R_n^*$ is.
Let $c_{d_1} = (2/\pi)^{1/2} + 0.59(d_1^{1/4} - 1)$. By Lemma
\ref{lem:gauss}, translating a convex Borel subset of $\R{d_1}$ by a vector
$v$ changes its standard normal probability by at most $c_{d_1}\|v\|$.

\begin{proposition} \label{prop:bound}
If Assumption \ref{assm:information} holds and $\Lambda(\psi_n)$ is closed and
convex, then
\[
  \sup_{C \in \mathcal{C}}\ \bigl| \pr(R^\ast_n \in C) - \pr(Z_\psi \in C) \bigr|
  \leq c_{d_1} \|\tilde\beta_n(\theta_n)\|\, \E\max(\|Z_\lambda\| - r_n, 0).
\]
\end{proposition}

With Corollary \ref{cor:naive}, the proposition bounds the asymptotic coverage
error of $C_n(\alpha)$, and the bound decays like a normal tail in $r_n$ because
$\|Z_\lambda\|^2 \sim \chi^2_{d_2}$. If $r_n \to \infty$ on $K$, the coverage of $C_n(\alpha)$ is asymptotically $1 - \alpha$ uniformly on
$K$ and in $\alpha$, as in Theorem \ref{thm:prof}.

We now consider the case $d_1 = d_2 = 1$ in detail. The local nuisance
set $\mathcal{H}_n$ is then an interval containing $0$, which we write as
\[
  \mathcal{H}_n = [-h_n, u_n], \qquad 0 \leq h_n, u_n \leq \infty,
\]
with an infinite endpoint excluded. Thus $h_n$ and $u_n$ are the standardized
distances from $\lambda_n$ to the lower and upper boundaries of
$\Lambda(\psi_n)$, $r_n = \min(h_n, u_n)$,
and $\Proj_n(z) = \min\{\max(z, -h_n), u_n\}$. Because
$\tilde{\beta}_n(\theta_n) = \kappa_n/(1 - \kappa_n^2)^{1/2}$ by Lemma
\ref{lem:Akappa},
\[
  R^\ast_n = Z_\psi + \kappa_n(1 - \kappa_n^2)^{-1/2}\{\min(Z_\lambda + h_n, 0)
  + \max(Z_\lambda - u_n, 0)\},
\]
whose distribution depends on $n$ only through $\kappa_n$, $h_n$ and $u_n$.
The distribution function is a one-dimensional integral of the standard
normal distribution function $\Phi$ against the normal density, so
probabilities and quantiles can be computed.

Suppose for now that $u_n = \infty$ and, to simplify notation, let $\kappa =
\kappa_n$ and $h = h_n$, so that
\[
  R^\ast_n = Z_\psi + \kappa(1 - \kappa^2)^{-1/2}\min(Z_\lambda + h, 0).
\]
Table~\ref{tab:cover} reports, for $\kappa < 0$ and several $h$, the tail
probabilities $\pr(R^\ast_n < -1.96)$ and $\pr(R^\ast_n > 1.96)$. By
Corollary \ref{cor:naive}, their sum approximates the non-coverage
probability of $C_n(0.05)$ at any $\theta$ with $\kappa_n(\theta) = \kappa$
and $h_n(\theta) = h$, with an error that vanishes uniformly in $\theta$ as
$n \to \infty$. Replacing $\psi$ by $-\psi$ changes the sign of $\kappa$ and
of $R^\ast_n$, so the sum is the same for $\kappa$ and $-\kappa$. In the table, the upper tail probability is above $0.025$ and the lower
below it for $h < 2$. The tails are asymmetric because, when $\kappa < 0$,
$R^\ast_n \geq Z_\psi$. At $h = 0$ the upper tail increases with $|\kappa|$
and is $0.28$ at $\kappa = -0.95$; the lower tail changes little with
$\kappa$. Both tails approach $0.025$ as $h$ increases, the more
slowly the larger $|\kappa|$, and the non-coverage probability, the sum of the
two, is within $0.005$ of $0.05$ for $h \geq 2$ at every $\kappa$ in the
table. The Size row is the sum at $h = 0$, the largest over $h$ in every
column.

\begin{table}[htbp]
\centering
\caption{
Probabilities $\pr(R^\ast_n < -1.96)$ (Lower) and $\pr(R^\ast_n > 1.96)$ (Upper). The size is the largest sum over the upper and lower probabilities over values of $h$.
}
\begin{tabular}{lrrrrrrrr}
	\toprule
	& \multicolumn{2}{c}{$\kappa = -0.45$} & \multicolumn{2}{c}{$\kappa = -0.70$} & \multicolumn{2}{c}{$\kappa = -0.85$} & \multicolumn{2}{c}{$\kappa = -0.95$} \\
	\cmidrule(lr){2-3} \cmidrule(lr){4-5} \cmidrule(lr){6-7} \cmidrule(lr){8-9}
	$h$ & Lower & Upper & Lower & Upper & Lower & Upper & Lower & Upper \\
	\midrule
	0.00 & 0.018 & 0.047 & 0.016 & 0.090 & 0.015 & 0.161 & 0.014 & 0.282 \\
	0.25 & 0.020 & 0.040 & 0.018 & 0.069 & 0.017 & 0.119 & 0.016 & 0.212 \\
	0.50 & 0.021 & 0.035 & 0.020 & 0.054 & 0.019 & 0.088 & 0.018 & 0.155 \\
	1.00 & 0.023 & 0.029 & 0.023 & 0.036 & 0.022 & 0.049 & 0.022 & 0.079 \\
	1.50 & 0.024 & 0.026 & 0.024 & 0.029 & 0.024 & 0.033 & 0.024 & 0.044 \\
	2.00 & 0.025 & 0.025 & 0.025 & 0.026 & 0.025 & 0.027 & 0.025 & 0.030 \\
	3.00 & 0.025 & 0.025 & 0.025 & 0.025 & 0.025 & 0.025 & 0.025 & 0.025 \\
	\midrule
	Size & \multicolumn{2}{c}{0.065} & \multicolumn{2}{c}{0.106} & \multicolumn{2}{c}{0.176} & \multicolumn{2}{c}{0.295} \\
	\bottomrule
\end{tabular}
\label{tab:cover}
\end{table}

Table~\ref{tab:cover} shows that normal critical values can give non-coverage
probabilities far above the nominal level near the boundary. One way to address this is by using least-favorable
critical values \citep{andrews2009hybrid}.  Consider the upper one-sided test of the null that the interest parameter
equals $\psi$, which rejects for large $\hat R_n(\psi)$. For $b, h \in \R{}$,
let $R^\ast_{b,h} = Z_\psi + b\min(Z_\lambda + h, 0)$, so that $R^\ast_n =
R^\ast_{b,h}$ with $b = \kappa(1 - \kappa^2)^{-1/2}$. For $\alpha \in (0, 1)$
and $\kappa \in (-1, 1)$, let $c_\alpha(\kappa)$ be the larger of
$z_{1-\alpha}$ and the $(1 - \alpha)$-quantile of $R^\ast_{b,0}$. When $b < 0$,
$R^\ast_{b,h}$ is at least $Z_\psi$ and stochastically decreasing in $h$; when
$b \geq 0$, it is at most $Z_\psi$. Hence $c_\alpha(\kappa)$ is the supremum
over $h \geq 0$ of the $(1 - \alpha)$-quantile of $R^\ast_{b,h}$, attained at
$h = 0$ when $\kappa < 0$. The least favorable critical value $c_\alpha\{\kappa_n(\hat\theta_n)\}$
replaces $\kappa_n$ by its estimator and takes the least favorable value
of $h$. The lower one-sided test rejects when $\hat{R}_n(\psi) <
-c_\alpha\{-\kappa_n(\hat\theta_n)\}$, obtained by applying the upper test to
$-\psi$. The next result shows that the test has asymptotic size at most
$\alpha$, and exactly $\alpha$ if $K$ contains a sequence approaching the
boundary along which $\kappa_n$ has a negative limit.

\begin{proposition}
\label{prop:lf}
If $d_1 = d_2 = 1$ and, on a set $K \subseteq \Theta$, Assumptions
\ref{assm:information}--\ref{assm:linear_approx},
\ref{assm:proj} and \ref{assm:bddcorr} hold, $u_n \to \infty$ and
$\kappa_n(\hat\theta_n) - \kappa_n = o_p(1)$, then
\[
  \limsup_{n \to \infty}\ \sup_{\theta \in K}\ \sup_{\alpha \in (0, 1)}\
  \bigl( \pr_\theta[\hat{R}_n(\psi) > c_\alpha\{\kappa_n(\hat\theta_n)\}] - \alpha \bigr) \leq 0 .
\]
If in addition $K$ contains a sequence along which $h_n \to 0$ and $\kappa_n$
converges to a negative limit, then $\limsup_{n \to \infty} \sup_{\theta \in K} \pr_\theta[\hat{R}_n(\psi) >
c_\alpha\{\kappa_n(\hat\theta_n)\}] = \alpha$, for every $\alpha \in (0, 1)$.
\end{proposition}

The proof is in Appendix \ref{sec:irreg_proofs}. It uses
$\kappa_n(\hat\theta_n)$ only through its consistency, so any consistent
estimator of $\kappa_n$ can replace it. Inverting the test gives
the confidence set
\[
  C_n^{\mathrm{lf}}(\alpha) = \{\psi : \hat{R}_n(\psi) \leq
  c_\alpha\{\kappa_n(\hat\theta_n)\}\},
\]
with $\hat\theta_n = \{\psi, \hat\lambda_n(\psi)\}$ as in \eqref{eq:eff_stat}.
It uses a partial maximizer over $\Lambda(\psi)$, with no extended set. Taking
complements in Proposition \ref{prop:lf}, its coverage is asymptotically at
least $1 - \alpha$, uniformly on $K$ and in $\alpha$.

To compare power, let $\theta_n = (\psi_n, \lambda_n)$ be the null value and
$(\psi_n + \Delta_n, \lambda_n)$ the true parameter, with $\Delta_n =
\bar{\mathcal{I}}_n(\theta_n)^{-1/2}\delta$ for a fixed $\delta \in \R{d_1}$.
To first order, with the expected information in place of the observed,
$S_n(\theta_n) \approx S_n(\psi_n + \Delta_n, \lambda_n) +
\mathcal{I}_n(\theta_n)(\Delta_n, 0)^\tsp$, whose second term adds $\delta$ to
the standardized efficient score and $\tilde\beta_n(\theta_n)^\tsp\delta$ to
the standardized nuisance score at $\theta_n$. Let $R^\ast_{n,\delta}$ be
\eqref{eq:Tstar} with $Z_\psi + \delta$ and $Z_\lambda +
\tilde\beta_n(\theta_n)^\tsp\delta$ in place of $Z_\psi$ and $Z_\lambda$. We
compare the tests under these approximating distributions.

If
$d_1 = d_2 = 1$ and $u_n = \infty$, then
\[
  R^\ast_{n,\delta} = \delta + Z_\psi + \kappa(1 - \kappa^2)^{-1/2}\min\{Z_\lambda
  + h + \kappa(1 - \kappa^2)^{-1/2}\delta, 0\},
\]
that is, $\delta$ plus $R^\ast_n$ with $h + \kappa(1 - \kappa^2)^{-1/2}\delta$
in place of $h$.
The power of the test of Proposition \ref{prop:lf} under the approximating
distribution is $\pr\{R^\ast_{n,\delta} > c_\alpha(\kappa)\}$, against $1 -
\Phi(z_{1-\alpha} - \delta)$ for the efficient score statistic of Theorem
\ref{thm:eff}, whose approximating distribution under the alternative is
$\rN(\delta, 1)$ for every $h$.

\begin{figure}[htbp]
\centering
\includegraphics[width=\textwidth]{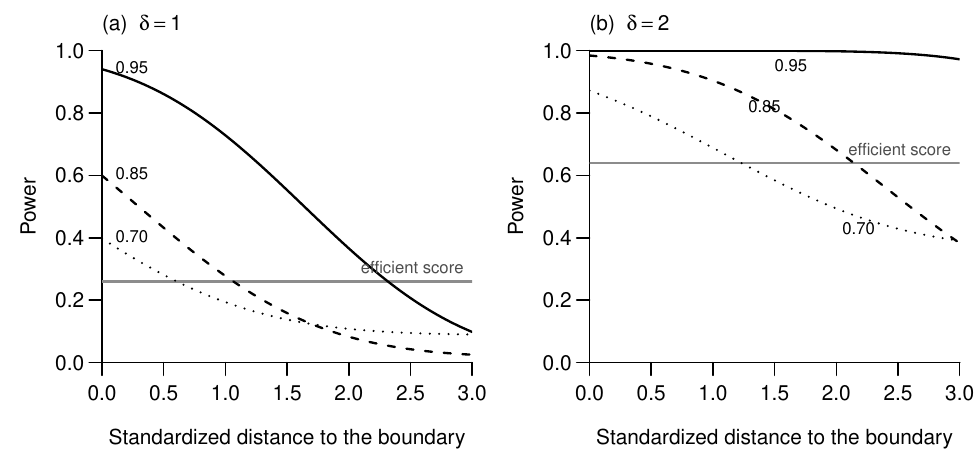}
\caption{Power of the one-sided test of Proposition \ref{prop:lf} at nominal level $0.05$, against the alternatives $\delta = 1$ and $\delta = 2$, as a function of the standardized distance to the boundary, computed as $\pr\{R^\ast_{n,\delta} > c_{0.05}(\kappa)\}$. The lines correspond to $\kappa = -0.70$, $-0.85$ and $-0.95$, labelled by $|\kappa|$, with critical values $c_{0.05}(\kappa)$ equal to $2.35$, $3.13$ and $5.27$, respectively. The gray line is the efficient score test, whose power does not depend on distance to the boundary.}
\label{fig:power}
\end{figure}

Figure \ref{fig:power} shows the two powers as functions of $h$ for
$\kappa = -0.70$, $-0.85$ and $-0.95$ and $\delta = 1, 2$. Applying the upper
test to $-\psi$ shows that, for $\kappa > 0$, the lower one-sided test
against $-\delta$ has the same power. In the figure, the test of Proposition
\ref{prop:lf} is more powerful than the efficient score test near the
boundary and less powerful beyond. At $h = 0$ its power against $\delta = 1$
is $0.40$, $0.60$ and $0.94$ for the three values of $\kappa$, against $0.26$
for the efficient score test, and the two powers cross at $h = 0.6$, $1.1$
and $2.3$. Against $\delta = 2$ the crossings are at $h = 1.2$ and $2.2$ for
the first two values of $\kappa$ and beyond the range of the figure for the
third.

The gain near the boundary is due to the alternative shifting the
nuisance score down when $\kappa < 0$ and $\delta > 0$, so that the
estimator is on the boundary more often than under the
null. On that event the statistic
shifts by $\delta/(1 - \kappa^2)$ rather than $\delta$. The factor is the
ratio of the full to the efficient information for $\psi$, so the gain grows
with $|\kappa|$. Far from the boundary the event is rare, the test is less
powerful because $c_\alpha(\kappa) > z_{1-\alpha}$, and its power tends to
$1 - \Phi\{c_\alpha(\kappa) - \delta\}$ as $h \to \infty$.

The lower one-sided test has critical value $-c_\alpha(-\kappa) =
-z_{1-\alpha}$, because $R^\ast_n \geq Z_\psi$ when $\kappa < 0$. The last
term of $R^\ast_{n,\delta}$ is nonnegative for every $\delta$, so
$R^\ast_{n,\delta} \geq \delta + Z_\psi$ and the lower test is never more
powerful than the efficient score test. In the configurations of the figure
it loses at most $0.03$. A two-sided test cannot gain in both directions. The
symmetric region $|\hat{R}_n(\psi)| > c$, with $c$ set under the null at
$h = 0$, is more powerful than the two-sided efficient score test near the
boundary when $\delta > 0$. Against $\delta = -1$ its power is at most $0.08$
in those configurations, against $0.17$. The gain is only for $\delta$ of
the opposite sign to $\kappa$ and for true parameters near the boundary.

\section{Examples} \label{sec:examples}

\subsection{Zero-inflated Poisson model} \label{sec:zip}

Let $Y_1, \dots, Y_n$ be independent and identically distributed from
\begin{equation}
\label{eq:zipmodel}
  \begin{aligned}
  \pr_\theta(Y = 0) &= \lambda + (1-\lambda)e^{-\psi}, \\
  \pr_\theta(Y = k) &= (1-\lambda)e^{-\psi}\psi^k/k!, \quad k \geq 1
  \end{aligned}
\end{equation}
where the parameter of interest $\psi$ is the Poisson mean and the nuisance parameter $\lambda$ is the zero-inflation probability, so $d_1 = d_2 = 1$. The parameter set is
\[
  \Theta = (0,\infty) \times [0,1), \qquad \Lambda(\psi) = [0,1), \quad \psi > 0,
\]
so that $\lambda = 0$, where the distribution reduces to $\mathrm{Poisson}(\psi)$, is a boundary point of $\Lambda(\psi)$.

The right-hand side of \eqref{eq:zipmodel} remains a probability mass function for small negative $\lambda$; with $\lambda$ of either sign it is the zero-modified Poisson model \citep{dietz2000estimation}. The nuisance estimator used to define the profile score statistic is
\begin{equation}
\label{eq:zipnuisdef}
  \hat\lambda_n(\psi) \in \argmax_{\lambda \in \tilde\Lambda(\psi)}
    \ell_n(\psi, \lambda),
  \qquad \hat\lambda_n = \hat\lambda_n(\psi_n),
\end{equation}
where $\tilde\Lambda$ is the extended set of Section \ref{sec:roots}, which
here does not depend on $n$,
\begin{equation}
\label{eq:zipext}
\begin{aligned}
  \tilde\Lambda(\psi) &= \bigl(-e^{-\psi}/\{1-e^{-\psi}\},\ 1\bigr), \\
  \tilde{\Theta} &= \bigl\{ \theta = (\psi,\lambda) : \psi > 0,\
    \lambda \in \tilde\Lambda(\psi) \bigr\}.
\end{aligned}
\end{equation}
The set $\tilde{\Theta}$ is open and contains $\Theta$, so $\lambda = 0$ is an
interior point of $\tilde\Lambda(\psi)$.

Quantities for one observation drop the subscript $n$. Write $S(y;\theta)$
for the score of an observation equal to $y$, so
\begin{equation}
\label{eq:zipsum}
  S_n(\theta) = \sum_{i=1}^n S(Y_i;\theta),
  \qquad
  \mathcal{I}_n(\theta) = n\,\mathcal{I}(\theta),
  \qquad
  \mathcal{I}(\theta) = \cov_\theta\{ S(Y;\theta) \},
\end{equation}
and $S$ and $\mathcal{I}$ have the same block notation as $S_n$ and
$\mathcal{I}_n$. The efficient information is likewise proportional to $n$,
\begin{equation}
\label{eq:zipeff}
  \bar{\mathcal{I}}_n(\theta) = n \bar{\mathcal{I}}(\theta)
  = n \left\{ \mathcal{I}_{\psi\psi}
    - \frac{\mathcal{I}_{\psi\lambda}^2}{\mathcal{I}_{\lambda\lambda}}
    \right\}(\theta) .
\end{equation}

The estimator \eqref{eq:zipnuisdef} takes negative values with probability
bounded away from zero along near-boundary sequences. Lemma
\ref{lem:zipinfo} shows that $\mathcal{I}$ is continuous and positive definite
on $\tilde{\Theta}$.

Let $n_0$ be the number of zeros in the sample. With probability tending
to one along every sequence in a compact compact set, $0 < n_0 < n$, and then $\lambda \mapsto
\ell_n(\psi,\lambda)$ is strictly concave on $\tilde\Lambda(\psi)$ with unique
maximizer $\hat\lambda_n(\psi) = (n_0/n - e^{-\psi})/(1 - e^{-\psi})$. Since $\tilde\Lambda(\psi_n)$ is open and
$\ell_n(\psi_n,\cdot)$ is differentiable on it, that maximizer is a root of the
nuisance score, so Assumption \ref{assm:nuis_est} holds. Lemma
\ref{lem:zipnuis} states this and the rate $\hat\lambda_n - \lambda_n =
O_p(n^{-1/2})$.

The distance $\tilde r_n$ of Section \ref{sec:roots} satisfies $\tilde r_n^2
\geq n\min\{\pr(Y = 0), \pr(Y > 0)\}$, as shown in the proof of the following
result, so $\tilde r_n \to \infty$ on every compact subset of $\Theta$. The proof
verifies the other hypotheses of Theorem \ref{thm:roots} and applies it.

\begin{theorem}
\label{thm:zip}
If $K \subseteq \Theta$ is compact and $\hat\lambda_n$ is as in
\eqref{eq:zipnuisdef}, then \eqref{eq:coverage} holds, and $\hat R_n(\psi_n)
\rightsquigarrow \rN(0, 1)$ on $K$.
\end{theorem}

The proof is in Appendix \ref{sec:zip_proofs}. We now take the nuisance
estimator to be a partial maximizer over $\Lambda(\psi)$,
\begin{equation}
\label{eq:zipconstr}
  \hat\lambda_n(\psi) \in \argmax_{\lambda \in \Lambda(\psi)} \ell_n(\psi,\lambda),
  \qquad \Lambda(\psi) = [0,1),
\end{equation}
in place of the partial maximizer \eqref{eq:zipnuisdef} over $\tilde\Lambda(\psi)$.
The standardized coefficient $\tilde\beta_n(\theta)$ of \eqref{eq:An} does
not depend on $n$, because the factors of $n$ in \eqref{eq:zipsum} and
\eqref{eq:zipeff} cancel, and it is finite by Lemma \ref{lem:zipinfo}; since
$d_1 = d_2 = 1$, $\tilde\beta_n(\theta) = \kappa_n(\theta)/\{1 -
\kappa_n(\theta)^2\}^{1/2}$. Let
\[
  h_n(\theta) = \mathcal{I}_{n\lambda\lambda}(\theta)^{1/2}\,\lambda = \lambda\bigl\{ n\, \mathcal{I}_{\lambda\lambda}(\theta) \bigr\}^{1/2}.
\]
By \eqref{eq:zipinfo} in Appendix \ref{sec:zip_proofs}, the squared
score correlation at $\lambda = 0$ is $\kappa_n(\theta)^2 = \psi/(e^{\psi} - 1)$, which decreases from one to zero as
$\psi$ increases, so that $\kappa_n(\theta) = -\{\psi/(e^{\psi} - 1)\}^{1/2}$
there: the smaller the Poisson mean, the harder it is to distinguish
zero-inflation from Poisson zeros, and the larger $|\kappa_n(\theta)|$. The value
$\kappa_n(\theta) = -0.85$ used below corresponds to $\psi = 0.62$ at $\lambda = 0$.

Assumption \ref{assm:proj} requires $\Lambda(\psi_n)$ to be closed and convex, and $[0,1)$ is
convex but not closed. We therefore fix $\bar\lambda \in (0,1)$ and take
$\Lambda(\psi) = [0,\bar\lambda]$ for every $\psi$ in \eqref{eq:zipconstr} for the rest of the subsection. The partial maximizer \eqref{eq:zipnuisdef} over $\tilde\Lambda(\psi_n)$ is less than one, as
the proof of Lemma \ref{lem:zipnuis} notes, so along every sequence with
$\limsup_n \lambda_n < \bar\lambda$ the change alters $\hat\lambda_n$ only on
an event of probability tending to zero.

\begin{lemma}
\label{lem:zipproj}
If $K \subseteq (0,\infty) \times [0, \bar\lambda)$ is compact and
$\hat\lambda_n = \hat\lambda_n(\psi_n)$ is as in \eqref{eq:zipconstr}, then Assumptions \ref{assm:information}--\ref{assm:linear_approx},
\ref{assm:coef_stable}, \ref{assm:proj} and \ref{assm:bddcorr} hold on $K$,
and the local nuisance set is the closed interval
\begin{equation}
\label{eq:zipH}
  \mathcal{H}_n = \bigl[-h_n(\theta_n),\ u_n\bigr],
  \qquad
  u_n = \bigl(\bar\lambda-\lambda_n\bigr)
        \mathcal{I}_{n\lambda\lambda}(\theta_n)^{1/2}
  \rightarrow \infty .
\end{equation}
\end{lemma}

By Lemma \ref{lem:zipproj}, Theorem \ref{thm:approx} applies on every compact
$K \subseteq (0, \infty) \times [0,\bar\lambda)$, with $\hat R_n$ the signed
root \eqref{eq:Rn}. Because $u_n$ in
\eqref{eq:zipH} is at least a constant times $n^{1/2}$ on $K$,
$\Proj_n(Z_\lambda) = \max(Z_\lambda, -h_n)$ with probability tending to
one, so the approximating variable is $R^\ast_n = Z_\psi + \kappa(1 -
\kappa^2)^{-1/2}\min(Z_\lambda + h, 0)$ with $\kappa = \kappa_n(\theta)$ and $h =
h_n(\theta)$, and its distribution depends on $\theta$ through $(\kappa, h)$ only. Along interior
sequences with $\lambda_n$ of order $n^{-1/2}$, $h_n$ stays bounded and
the approximating distribution has nonzero mean $\kappa_n(\theta_n)\{1 - \kappa_n(\theta_n)^2\}^{-1/2}\E\{\min(Z_\lambda + h_n, 0)\}$, so the pointwise normal limit that holds at each fixed interior $\theta$ fails along such sequences.

Theorem \ref{thm:eff} also applies with the estimator \eqref{eq:zipconstr},
by Lemma \ref{lem:zipproj}.

\begin{figure}[htbp]
\centering
\includegraphics[width=\textwidth]{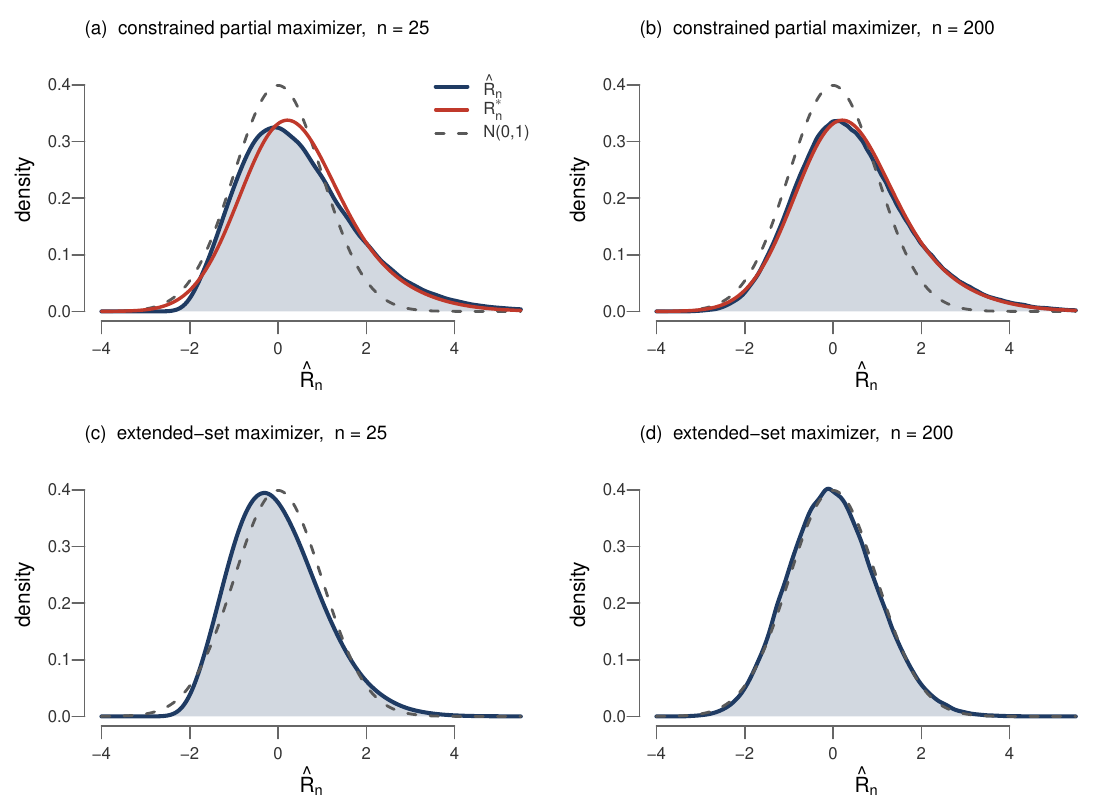}
\caption{Sampling density of $\hat R_n$ ($200{,}000$ replicates, blue and
shaded), the density of $R^\ast_n$ with $\kappa = -0.85$ and $h = 0.25$ (red),
and the $\rN(0,1)$ density (dashed). The top row uses the partial
maximizer \eqref{eq:zipconstr} over $\Lambda(\psi)$ and the bottom row the partial maximizer \eqref{eq:zipnuisdef} over $\tilde\Lambda(\psi)$. The columns are
$n = 25$ and $n = 200$, with $(\psi, \lambda) = (0.5965, 0.05501)$ and
$(0.6103, 0.01924)$, chosen so that $\kappa_n(\theta) = -0.85$ and
$h_n(\theta) = 0.25$ in both.}
\label{fig:zipkde}
\end{figure}

Figure \ref{fig:zipkde} fixes $\kappa_n(\theta) = -0.85$ and $h_n(\theta) = 0.25$ at
$n = 25$ and $n = 200$. The distributions of $R^\ast_n$ and of the standard normal
are then the same in the two columns, so any difference between
the columns is due to $n$. The top row uses the partial maximizer
\eqref{eq:zipconstr} over $\Lambda(\psi)$ and the bottom row the partial maximizer \eqref{eq:zipnuisdef} over $\tilde\Lambda(\psi)$, with which the
statistic is asymptotically standard normal by Theorem \ref{thm:zip}. In the top row the density of $\hat R_n$ is right-skewed, with a lower peak and a heavier right tail than the standard normal at both sample sizes, and the density of $R^\ast_n$ follows it, closely at $n = 200$ and with a slightly shifted peak at $n = 25$. In the bottom row the density is close to the standard normal at $n = 200$ and shifted slightly to the left at $n = 25$.

Because $\kappa_n(\theta) < 0$ in this model, the gain in Figure
\ref{fig:power} is for alternatives with the true Poisson mean above the
hypothesized one, that is, for the lower confidence bound on the mean. The
dashed lines are at the value of $\kappa$ this model has for $\psi$ near
$0.6$. The consistency
$\kappa_n(\hat\theta_n) - \kappa_n(\theta_n) = o_p(1)$ that Proposition
\ref{prop:lf} requires follows from continuity of $\kappa_n$ (equivalently, of $\tilde\beta_n$)
and the rate in the proof of Lemma \ref{lem:zipproj}.

Table~\ref{tab:zipcov} gives the coverage of nominal $95\%$ regions for
$\psi$ at parameter values with the zero-inflation probability at or near its
boundary. The region of Theorem~\ref{thm:zip}, which uses the partial maximizer
over $\tilde\Lambda(\psi)$, and the efficient score region of
Theorem~\ref{thm:eff}, which uses the partial maximizer over $\Lambda(\psi)$, are within
\zipvalidmaxdev{} of nominal in every cell. The naive region, the chi-squared
region with the partial maximizer over $\Lambda(\psi)$, has coverage as low as \zipnaivemin, and in every cell its coverage is within \zippreddiff{} of the
asymptotic value of Corollary~\ref{cor:naive}. Its lower confidence bound has coverage as low as \zipnaiveonemin, and in every cell within \zippredonediff{} of its asymptotic value. Wald and likelihood ratio intervals with the partial maximizer over $\Lambda(\psi)$
and the same critical values have coverage between \zipwaldmin{} and
\zipwaldmax{} and between \ziplrtmin{} and \ziplrtmax{} in the same
configurations. The bound
$C_n^{\mathrm{lf}}(\alpha)$ of Proposition~\ref{prop:lf} has coverage at least
\ziplfmin, near nominal at the boundary and conservative away
from it.

\begin{table}[htbp]
\centering
\caption{
Coverage of nominal $95\%$ regions for the Poisson mean $\psi$ in the
zero-inflated Poisson model (\zipcovnsim{} samples per row). The columns Extended (Theorem~\ref{thm:prof}), Efficient (Theorem~\ref{thm:eff}),
Naive (partial maximizer over $\Lambda(\psi)$, normal critical values) and LF (the least favorable
bound of Proposition~\ref{prop:lf}), are simulated coverage probabilities.
Approx is the asymptotic coverage of the naive region calculated from the approximating sequence of distributions $R_n^*$
(Corollary~\ref{cor:naive}).
}
{\footnotesize\setlength{\tabcolsep}{4pt} 
	\begin{tabular}{rrrrrrrrrrr}
		\toprule
		& & & & \multicolumn{4}{c}{Two-sided} & \multicolumn{3}{c}{Lower bound} \\
		\cmidrule(lr){5-8} \cmidrule(lr){9-11}
		$n$ & $\lambda$ & $\tilde\beta_n(\theta)$ & $h_n(\theta)$ & Extended & Efficient & Naive & Approx & Naive & Approx & LF \\
		\midrule
		\multicolumn{11}{l}{$\psi = 0.5$} \\
		25 & 0.00 & $-1.83$ & 0.00 & 0.963 & 0.954 & 0.792 & 0.801 & 0.782 & 0.764 & 0.946 \\
		100 & 0.00 & $-1.83$ & 0.00 & 0.953 & 0.948 & 0.803 & 0.801 & 0.763 & 0.764 & 0.943 \\
		400 & 0.00 & $-1.83$ & 0.00 & 0.949 & 0.946 & 0.806 & 0.801 & 0.770 & 0.764 & 0.949 \\
		25 & 0.02 & $-1.82$ & 0.08 & 0.962 & 0.955 & 0.808 & 0.818 & 0.797 & 0.784 & 0.952 \\
		100 & 0.02 & $-1.82$ & 0.16 & 0.952 & 0.948 & 0.834 & 0.833 & 0.797 & 0.801 & 0.957 \\
		400 & 0.02 & $-1.82$ & 0.32 & 0.951 & 0.949 & 0.864 & 0.859 & 0.837 & 0.832 & 0.973 \\
		25 & 0.05 & $-1.80$ & 0.20 & 0.962 & 0.956 & 0.832 & 0.841 & 0.820 & 0.811 & 0.961 \\
		100 & 0.05 & $-1.80$ & 0.41 & 0.953 & 0.950 & 0.870 & 0.872 & 0.841 & 0.847 & 0.973 \\
		400 & 0.05 & $-1.80$ & 0.81 & 0.950 & 0.949 & 0.912 & 0.913 & 0.896 & 0.898 & 0.990 \\
		25 & 0.20 & $-1.73$ & 0.85 & 0.960 & 0.959 & 0.915 & 0.917 & 0.897 & 0.904 & 0.986 \\
		100 & 0.20 & $-1.73$ & 1.69 & 0.953 & 0.953 & 0.948 & 0.945 & 0.932 & 0.943 & 0.996 \\
		400 & 0.20 & $-1.73$ & 3.39 & 0.950 & 0.950 & 0.950 & 0.950 & 0.946 & 0.950 & 0.999 \\
		\addlinespace
		\multicolumn{11}{l}{$\psi = 1.25$} \\
		25 & 0.00 & $-1.00$ & 0.00 & 0.956 & 0.949 & 0.882 & 0.892 & 0.851 & 0.860 & 0.940 \\
		100 & 0.00 & $-1.00$ & 0.00 & 0.951 & 0.948 & 0.894 & 0.892 & 0.852 & 0.860 & 0.946 \\
		400 & 0.00 & $-1.00$ & 0.00 & 0.950 & 0.949 & 0.890 & 0.892 & 0.856 & 0.860 & 0.949 \\
		25 & 0.02 & $-0.98$ & 0.16 & 0.955 & 0.949 & 0.893 & 0.907 & 0.868 & 0.880 & 0.947 \\
		100 & 0.02 & $-0.98$ & 0.31 & 0.949 & 0.947 & 0.913 & 0.917 & 0.882 & 0.895 & 0.961 \\
		400 & 0.02 & $-0.98$ & 0.62 & 0.949 & 0.949 & 0.927 & 0.931 & 0.911 & 0.918 & 0.974 \\
		25 & 0.05 & $-0.95$ & 0.38 & 0.956 & 0.952 & 0.913 & 0.922 & 0.893 & 0.903 & 0.960 \\
		100 & 0.05 & $-0.95$ & 0.76 & 0.951 & 0.950 & 0.935 & 0.937 & 0.916 & 0.926 & 0.976 \\
		400 & 0.05 & $-0.95$ & 1.53 & 0.951 & 0.951 & 0.948 & 0.948 & 0.942 & 0.945 & 0.986 \\
		25 & 0.20 & $-0.82$ & 1.44 & 0.955 & 0.954 & 0.949 & 0.948 & 0.939 & 0.945 & 0.976 \\
		100 & 0.20 & $-0.82$ & 2.88 & 0.951 & 0.951 & 0.951 & 0.950 & 0.946 & 0.950 & 0.982 \\
		400 & 0.20 & $-0.82$ & 5.77 & 0.949 & 0.949 & 0.949 & 0.950 & 0.948 & 0.950 & 0.983 \\
		\bottomrule
	\end{tabular}
}
\label{tab:zipcov}
\end{table}

Table~\ref{tab:zipbox} gives the simultaneous coverage of boxes for $(\psi,
\lambda)$ from the same samples. The interval for $\lambda$ reverses the roles of the
two parameters: it uses a root of the score for $\psi$ over $\Psi(\lambda)$, at
the hypothesized $\lambda$, as the nuisance estimator. The box of
Proposition~\ref{prop:box} uses $\Omega_n$ at the maximum likelihood estimator
over $\Theta$. The interval for
$\lambda$, whose true value is on or near the boundary of its set, has coverage
between \zippimin{} and \zippimax. The product of the two marginal intervals
has joint coverage as low as \zipboxpmin, the Bonferroni box between
\zipboxbmin{} and \zipboxbmax, and the box of Proposition~\ref{prop:box}
between \zipboxmmin{} and \zipboxmmax, with average critical value between
\zipmbarmin{} and \zipmbarmax{} against $1.96$ for the product and $2.24$ for the Bonferroni box.

\begin{table}[htbp]
\centering
\caption{
Simultaneous coverage of nominal $95\%$ boxes for $(\psi, \lambda)$ in the
zero-inflated Poisson model (\zipcovnsim{} samples per row). ``Marginal'' is the
coverage of the chi-squared interval for $\lambda$, with $\psi$ estimated by a root
of its score. The three boxes combine that interval with the Extended interval for
$\psi$ of Table~\ref{tab:zipcov}, using critical value $1.96$ (Product), $2.24$
(Bonferroni) or $m_{0.05}(\hat\Omega_n)$ (Joint, Proposition~\ref{prop:box}), which
varies with the sample, and we report the average $\bar m$.
}
{\footnotesize\setlength{\tabcolsep}{4pt}
	\begin{tabular}{rrrrrrrr}
		\toprule
		& & & & \multicolumn{1}{c}{Marginal} & \multicolumn{3}{c}{Box} \\
		\cmidrule(lr){5-5} \cmidrule(lr){6-8}
		$n$ & $\lambda$ & $|\kappa|$ & $\bar m$ & $\lambda$ & Product & Bonferroni & Joint \\
		\midrule
		\multicolumn{8}{l}{$\psi = 0.5$} \\
		25 & 0.00 & 0.88 & 2.13 & 0.963 & 0.941 & 0.966 & 0.951 \\
		100 & 0.00 & 0.88 & 2.13 & 0.950 & 0.927 & 0.963 & 0.952 \\
		400 & 0.00 & 0.88 & 2.12 & 0.950 & 0.928 & 0.962 & 0.952 \\
		25 & 0.02 & 0.88 & 2.13 & 0.953 & 0.930 & 0.966 & 0.954 \\
		100 & 0.02 & 0.88 & 2.13 & 0.949 & 0.927 & 0.964 & 0.953 \\
		400 & 0.02 & 0.88 & 2.12 & 0.951 & 0.929 & 0.964 & 0.953 \\
		25 & 0.05 & 0.87 & 2.13 & 0.953 & 0.930 & 0.963 & 0.956 \\
		100 & 0.05 & 0.87 & 2.13 & 0.950 & 0.928 & 0.963 & 0.952 \\
		400 & 0.05 & 0.87 & 2.12 & 0.951 & 0.927 & 0.962 & 0.951 \\
		25 & 0.20 & 0.87 & 2.13 & 0.956 & 0.931 & 0.962 & 0.953 \\
		100 & 0.20 & 0.87 & 2.13 & 0.950 & 0.927 & 0.960 & 0.951 \\
		400 & 0.20 & 0.87 & 2.13 & 0.951 & 0.927 & 0.962 & 0.950 \\
		\addlinespace
		\multicolumn{8}{l}{$\psi = 1.25$} \\
		25 & 0.00 & 0.71 & 2.19 & 0.946 & 0.917 & 0.960 & 0.954 \\
		100 & 0.00 & 0.71 & 2.18 & 0.950 & 0.918 & 0.959 & 0.952 \\
		400 & 0.00 & 0.71 & 2.18 & 0.949 & 0.917 & 0.958 & 0.951 \\
		25 & 0.02 & 0.70 & 2.19 & 0.953 & 0.924 & 0.958 & 0.956 \\
		100 & 0.02 & 0.70 & 2.18 & 0.948 & 0.915 & 0.956 & 0.950 \\
		400 & 0.02 & 0.70 & 2.18 & 0.949 & 0.915 & 0.957 & 0.950 \\
		25 & 0.05 & 0.69 & 2.19 & 0.946 & 0.917 & 0.959 & 0.955 \\
		100 & 0.05 & 0.69 & 2.18 & 0.949 & 0.916 & 0.956 & 0.951 \\
		400 & 0.05 & 0.69 & 2.18 & 0.950 & 0.917 & 0.957 & 0.950 \\
		25 & 0.20 & 0.63 & 2.19 & 0.951 & 0.918 & 0.957 & 0.950 \\
		100 & 0.20 & 0.63 & 2.19 & 0.950 & 0.915 & 0.955 & 0.948 \\
		400 & 0.20 & 0.63 & 2.19 & 0.949 & 0.913 & 0.955 & 0.950 \\
		\bottomrule
	\end{tabular}
}
\label{tab:zipbox}
\end{table}

\subsection{One-way random effects model} \label{sec:oneway}

\citet{shedden2026score} verify Assumptions
\ref{assm:information}--\ref{assm:nuis_est} in linear mixed models under a
condition that holds for independent clusters and for crossed random effects.
However, they do not consider partial maximizers over $\Lambda(\psi)$, so the
approximating distribution and the least favorable test of Section
\ref{sec:irreg}, together with Assumptions \ref{assm:proj} and
\ref{assm:bddcorr}, are not treated there. The one-way model is simple enough that all of these
are in closed form, so everything in Section \ref{sec:irreg} can be computed
exactly.

Let $Y_{ij}$ be the $j$th measurement on unit $i$, for $i \in \{1, \dots, m\}$
and $j \in \{1, \dots, k\}$, with
\[
  Y_{ij} = W_i + E_{ij},
\]
where the unit effects $W_i \sim \rN(0,
\lambda)$ and errors $E_{ij} \sim \rN(0, \psi)$ are all independent. Two
measurements on the same unit have covariance $\lambda$, and different units
are independent. Here $k \geq 2$ is fixed, $n = mk$, and asymptotics are as
$m \to \infty$. The parameter of interest is $\psi > 0$, the
nuisance parameter is $\lambda \geq 0$, and $\Theta = (0,
\infty) \times [0, \infty)$, so $\Lambda(\psi) = [0, \infty)$ for every $\psi$.
Precision studies with $k$ replicates on each of $m$ units use this design; they
report an upper confidence limit for $\psi$ \citep{burdick1992confidence}, and the between-unit variance is
often estimated to be zero. 

Inference is based on the restricted likelihood, which does not depend on a common mean if one is added. The restricted likelihood depends on the data only through
\[
  Q_{\mathrm{b}} = k\sum_{i=1}^m (\bar Y_{i\cdot} - \bar Y_{\cdot\cdot})^2, \qquad
  Q_{\mathrm{w}} = \sum_{i=1}^m\sum_{j=1}^k (Y_{ij} - \bar Y_{i\cdot})^2,
\]
which are independent with $Q_{\mathrm{b}} \sim \tau\chi^2_{m-1}$ and $Q_{\mathrm{w}} \sim
\psi\chi^2_{m(k-1)}$, where $\tau = \tau(\theta) = \psi + k\lambda$. Up to a
constant, the restricted log-likelihood is
\begin{equation} \label{eq:owll}
  \ell_n(\theta) = -\tfrac12\{(m-1)\log\tau + Q_{\mathrm{b}}/\tau + m(k-1)\log\psi +
  Q_{\mathrm{w}}/\psi\},
\end{equation}
which is defined on the extended set $\tilde\Theta = \{\theta : \psi > 0,\
\psi + k\lambda > 0\}$, the set on which the covariance matrix of $Y_i = (Y_{i1}, \dots, Y_{ik})$ is positive definite, so $\tilde\Lambda(\psi) =
(-\psi/k, \infty)$. Differentiating \eqref{eq:owll},
\begin{equation} \label{eq:owscore}
  S_{n\psi}(\theta) = \frac{Q_{\mathrm{b}} - (m-1)\tau}{2\tau^2} + \frac{Q_{\mathrm{w}} -
  m(k-1)\psi}{2\psi^2}, \qquad
  S_{n\lambda}(\theta) = \frac{k\{Q_{\mathrm{b}} - (m-1)\tau\}}{2\tau^2},
\end{equation}
and, since $\var(Q_{\mathrm{b}}) = 2(m-1)\tau^2$ and $\var(Q_{\mathrm{w}}) = 2m(k-1)\psi^2$,
\begin{equation} \label{eq:owinfo}
  \mathcal{I}_{n\psi\psi} = \frac{m-1}{2\tau^2} + \frac{m(k-1)}{2\psi^2}, \qquad
  \mathcal{I}_{n\psi\lambda} = \frac{k(m-1)}{2\tau^2}, \qquad
  \mathcal{I}_{n\lambda\lambda} = \frac{k^2(m-1)}{2\tau^2}.
\end{equation}
By \eqref{eq:owinfo}, the efficient information is
$\bar{\mathcal{I}}_n = m(k-1)/(2\psi^2)$, the information in $Q_{\mathrm{w}}$
alone, and the regression coefficient $\beta_n(\theta) = 1/k$ does not depend on
$\theta$, so the efficient score is $U_n(\theta) = \{Q_{\mathrm{w}} -
m(k-1)\psi\}/(2\psi^2)$ for every value of the nuisance parameter. The
squared score correlation is
\[
  \kappa_n(\theta)^2 = \frac{(m-1)\psi^2}{(m-1)\psi^2 + m(k-1)\tau^2},
\]
which decreases in $\lambda$ from $(m-1)/(mk-1)$ at the boundary and never exceeds
$1/k$, so Assumption \ref{assm:bddcorr} holds on $\Theta$. By
the identity before Assumption \ref{assm:bddcorr},
\begin{equation} \label{eq:owa}
  \tilde\beta_n(\theta) = \frac{\kappa_n(\theta)}{\{1 - \kappa_n(\theta)^2\}^{1/2}} = \frac{\psi}{\tau}\Bigl\{\frac{m-1}{m(k-1)}\Bigr\}^{1/2}
  > 0,
\end{equation}
which at the boundary tends to $(k-1)^{-1/2}$. With $k = 2$, the
score correlation at the boundary tends to $2^{-1/2} \approx 0.71$, corresponding, up to rounding, to the column
$\kappa = -0.70$ of Table \ref{tab:cover}, with its two tails swapped because $\kappa > 0$ here. After reflection, the dotted lines of Figure \ref{fig:power} are at the value of $\kappa$ this model has at the boundary.

Fix $\psi$. Since $Q_{\mathrm{b}} > 0$, the
function $\tau \mapsto (m-1)\log\tau + Q_{\mathrm{b}}/\tau$ decreases
on $(0, Q_{\mathrm{b}}/(m-1)]$ and increases after, so the partial maximizer over $\tilde\Lambda(\psi)$ is $\{Q_{\mathrm{b}}/(m-1) - \psi\}/k$. The partial maximizer over $\Lambda(\psi)$ is
\begin{equation} \label{eq:owconstr}
  \hat\lambda_n(\psi) = \max\bigl[\{Q_{\mathrm{b}}/(m-1) - \psi\}/k,\ 0\bigr].
\end{equation}
Because $\tilde\Lambda(\psi)$ is open, the partial maximizer over it is a root of the
nuisance score, so Assumption \ref{assm:nuis_est} holds. Neither $U_n$ nor $\bar{\mathcal{I}}_n$ depends on
$\lambda$, and $\bar{\mathcal{I}}_n$ and $\mathcal{I}_{n\lambda\lambda}$ are positive on $\tilde\Theta$ by \eqref{eq:owinfo}, so Assumptions \ref{assm:eff_info_stable} and
\ref{assm:linear_approx} hold for every estimator with $\hat\lambda_n(\psi) \in \tilde\Lambda(\psi)$, and since
$\beta_n$ does not depend on $\theta$, the left-hand side of Assumption
\ref{assm:coef_stable} is zero, so Theorem \ref{thm:eff} applies with any
such estimator. The proof of Proposition
\ref{prop:oneway} verifies Assumptions \ref{assm:information} and
\ref{assm:clt} on $\Theta$. The boundary of
$\tilde\Lambda(\psi_n)$ is at distance $\tilde r_n =
\mathcal{I}_{n\lambda\lambda}(\theta_n)^{1/2}(\lambda_n + \psi_n/k) = \{(m-1)/2\}^{1/2}$
from every parameter value on the information scale.
Theorem \ref{thm:prof} therefore gives confidence regions for $\psi$,
formed with the partial maximizer over $\tilde\Lambda(\psi)$, with asymptotically correct
coverage probability uniformly on $\Theta$. For the partial maximizer \eqref{eq:owconstr} over $\Lambda(\psi)$, the standardized
distance to the boundary of Section \ref{sec:irreg} is $h_n(\theta)
= \mathcal{I}_{n\lambda\lambda}(\theta)^{1/2}\lambda = k\lambda\{(m-1)/2\}^{1/2}/\tau$.

\begin{proposition} \label{prop:oneway}
If $K \subseteq \Theta$ is compact and $\hat\lambda_n = \hat\lambda_n(\psi_n)$
is the partial maximizer \eqref{eq:owconstr} over $\Lambda(\psi)$, then Assumptions
\ref{assm:information}--\ref{assm:linear_approx}, \ref{assm:proj} and
\ref{assm:bddcorr} hold on $K$,
with $\tilde\beta_n$ as in \eqref{eq:owa} and $\mathcal{H}_n = [-h_n, \infty)$, and
$\tilde\beta_n(\hat\theta_n) - \tilde\beta_n(\theta_n) = o_p(1)$ on $K$.
\end{proposition}

The proof is in Appendix \ref{sec:oneway_proofs}. By Proposition
\ref{prop:oneway}, Theorem \ref{thm:approx} applies on every compact $K
\subseteq \Theta$. Because $u_n = \infty$, $R^\ast_n = Z_\psi + \kappa(1 -
\kappa^2)^{-1/2}\min(Z_\lambda + h, 0)$ with $\kappa = \kappa_n(\theta)$ and $h = h_n(\theta)$, and Proposition
\ref{prop:lf} applies.
Because $\kappa_n(\theta) > 0$, the partial maximizer over $\Lambda(\psi)$ gives more power when the true error
variance is below the hypothesized one, that is, for the upper confidence limit. On the event that
\eqref{eq:owconstr} is zero, the profile score is $\{Q_{\mathrm{b}} + Q_{\mathrm{w}} -
(mk-1)\psi\}/(2\psi^2)$, the score for $\psi$ when $\lambda = 0$ is known,
with information $(mk-1)/(2\psi^2)$ in place of
$m(k-1)/(2\psi^2)$. The ratio is $1/\{1 - \kappa_n(\theta)^2\}$, which tends to $k/(k-1)$. For $k=2$, the profile score on that event has about twice the efficient
information about the error variance.

\bibliographystyle{apalike}
\bibliography{profile_score_bdry}

\appendix
\counterwithin{theorem}{section}
\counterwithin{lemma}{section}
\counterwithin{corollary}{section}
\renewcommand{\theHtheorem}{\thesection.\arabic{theorem}}
\renewcommand{\theHlemma}{\thesection.\arabic{lemma}}
\renewcommand{\theHcorollary}{\thesection.\arabic{corollary}}
\counterwithin{equation}{section}
\renewcommand{\theHequation}{\thesection.\arabic{equation}}

\section{Proofs} \label{sec:proofs}

Statements of results write arguments as in the main text. In the proofs, a
score, information matrix or quantity defined from them in Sections
\ref{sec:introduction} and \ref{sec:model} that is written without an argument is evaluated at
$\theta_n$.


\subsection{Proofs for Section \ref{sec:general_theory}}
\label{sec:chi2_proofs}

\begin{proof}[Proof of Lemma~\ref{lem:unif}]
Both conclusions hold if $K$ is empty, so suppose it is not. If the
conclusion of (i) fails, there is an $\eta > 0$ such that
$\sup_{\theta \in K}|g_n(\theta) - c| > \eta$ for infinitely many $n$. For
those $n$ choose $\theta_n \in K$ with $|g_n(\theta_n) - c| > \eta$, and for
the other $n$ choose $\theta_n \in K$ arbitrarily. Then $(\theta_n)$ is a
sequence in $K$ along which $g_n(\theta_n)$ does not converge to $c$,
contradicting the hypothesis. If the conclusion of (ii) fails, there is a
$c' > c$ such that $\sup_{\theta \in K} g_n(\theta) > c'$ for infinitely many
$n$. Choosing $\theta_n$ as before, with $g_n(\theta_n) > c'$ for those $n$,
gives $\limsup_{n \to \infty} g_n(\theta_n) \geq c' > c$, contradicting the
hypothesis. This completes the proof.
\end{proof}

\begin{lemma}\label{lem:sqrt}
If $P$ is a symmetric positive definite $d_1 \times d_1$ matrix, $Q$ is a
symmetric invertible $d_1 \times d_1$ matrix, and
$\|P^{1/2}Q^{-1}P^{1/2} - I_{d_1}\| \leq \delta$ for some $\delta < 1$, then
$Q$ is positive definite and
$\|Q^{-1/2}P^{1/2} - I_{d_1}\| \leq (3 d_1 \delta)^{1/2}$.
\end{lemma}

\begin{proof}
The matrix $V = P^{1/2}Q^{-1}P^{1/2}$ is symmetric with eigenvalues in
$[1 - \delta, 1 + \delta]$, so it is positive definite, and hence so are
$Q^{-1} = P^{-1/2}VP^{-1/2}$ and $Q$. Let $M = Q^{-1/2}P^{1/2}$, so that
$M^\tsp M = V$ and the smallest singular
value of $M$, the square root of the smallest eigenvalue of $V$, is at least
$(1 - \delta)^{1/2}$. Because $P^{1/4}MP^{-1/4} = P^{1/4}Q^{-1/2}P^{1/4}$ is
symmetric positive definite, the eigenvalues of $M$ are positive, and each is
at least the smallest singular value of $M$
\citep[Theorem 5.6.9(b)]{horn2012matrix}. Hence
$\tr(M) \geq d_1(1 - \delta)^{1/2} \geq d_1(1 - \delta)$, and since
$\tr(M^\tsp M) = \tr(V) \leq d_1(1 + \delta)$,
\[
  \|M - I_{d_1}\|^2 \leq \tr\{(M - I_{d_1})^\tsp (M - I_{d_1})\}
  = \tr(V) - 2\tr(M) + d_1 \leq 3 d_1 \delta,
\]
where the first inequality holds because the spectral norm is at most the
Frobenius norm. This completes the proof.
\end{proof}

\begin{proof}[Proof of Theorem \ref{thm:roots}]
For every $M > 0$, because $\tilde r_n \to \infty$, $B_n^M = \{\lambda \in
\R{d_2} : \|\mathcal{I}_{n\lambda\lambda}^{1/2}(\lambda - \lambda_n)\| \leq M\}$
is in the interior of $\tilde\Lambda_n(\psi_n)$ for all large $n$, so
$S_n(\psi_n, \cdot)$ is continuously differentiable on $B_n^M$. Let $\veps_n^M$
and $\delta_n^M$ be the two suprema in \eqref{eq:info_unif}, so that $\veps_n^M
= o_p(1)$ and $\delta_n^M \to 0$.

To prove part (i), let $\hat\lambda_n \in \tilde\Lambda_n(\psi_n)$ satisfy
\eqref{eq:info_bounded}. It suffices to show that, for every fixed $M > 0$, on
$\{\hat\lambda_n \in B_n^M\}$ and for all large $n$, the quantities in
Assumptions \ref{assm:eff_info_stable} and \ref{assm:linear_approx} exist and
have norms at most $\delta_n^M/(1 - \delta_n^M)$ and $M\veps_n^M$. Indeed, in
that case, for any $\eta \in (0, 1)$, \eqref{eq:info_bounded} gives an $M$ with
$\pr(\hat\lambda_n \in B_n^M) \geq 1 - \eta$ for all large $n$, on which event
the quantities exist and are bounded by $\delta_n^M/(1 - \delta_n^M) \to 0$ and
$M\veps_n^M = o_p(1)$; since $\eta$ is arbitrary, both quantities are $o_p(1)$.

The matrix $\Gamma_n =
\bar{\mathcal{I}}_n^{-1/2}[I_{d_1}, -\beta_n]\mathcal{I}_n^{1/2}$ of the proof of
Lemma \ref{lem:chi2} satisfies $\Gamma_n\Gamma_n^\tsp = I_{d_1}$, and $L_n =
\bar{\mathcal{I}}_n^{1/2}[I_{d_1}, 0]\mathcal{I}_n^{-1/2}$ satisfies $L_nL_n^\tsp =
\bar{\mathcal{I}}_n^{1/2}\bar{\mathcal{I}}_n^{-1}\bar{\mathcal{I}}_n^{1/2} = I_{d_1}$
because $\bar{\mathcal{I}}_n^{-1}$ is the leading block of $\mathcal{I}_n^{-1}$;
hence both have spectral norm one. On $\{\hat\lambda_n \in B_n^M\}$ the
segment from $\lambda_n$ to $\hat\lambda_n$ is in $B_n^M$, and since
$\partial U_n(\psi_n, \lambda)/\partial\lambda^\tsp = -[I_{d_1},
-\beta_n]\mathcal{J}_{n\cdot\lambda}(\psi_n, \lambda)$ by \eqref{eq:eff_score},
\[
  U_n(\psi_n, \hat\lambda_n) - U_n(\theta_n) = -[I_{d_1}, -\beta_n]\int_0^1
  \mathcal{J}_{n\cdot\lambda}(\psi_n, \lambda_n(t))\,\dd t\; (\hat\lambda_n - \lambda_n),
\]
where $\lambda_n(t) = \lambda_n + t(\hat\lambda_n - \lambda_n)$. The matrix
$\mathcal{I}_{n\cdot\lambda}$ can be subtracted from the integrand because
$[I_{d_1}, -\beta_n]\mathcal{I}_{n\cdot\lambda} = \mathcal{I}_{n\psi\lambda} -
\beta_n\mathcal{I}_{n\lambda\lambda} = 0$. Multiplying by
$\bar{\mathcal{I}}_n^{-1/2}$ and using $\bar{\mathcal{I}}_n^{-1/2}[I_{d_1},
-\beta_n] = \Gamma_n\mathcal{I}_n^{-1/2}$ gives
\begin{align*}
  \tilde U_n(\psi_n, \hat\lambda_n) - \tilde U_n(\theta_n)
  &= -\Gamma_n\int_0^1 \mathcal{I}_n^{-1/2}\bigl\{\mathcal{J}_{n\cdot\lambda}(\psi_n, \lambda_n(t)) -
  \mathcal{I}_{n\cdot\lambda}\bigr\}\mathcal{I}_{n\lambda\lambda}^{-1/2}\,\dd t \\
  &\qquad \cdot \mathcal{I}_{n\lambda\lambda}^{1/2}(\hat\lambda_n - \lambda_n),
\end{align*}
whose norm is at most
$\veps_n^M\|\mathcal{I}_{n\lambda\lambda}^{1/2}(\hat\lambda_n -
\lambda_n)\| \leq M\veps_n^M$. For Assumption \ref{assm:eff_info_stable},
let $\Sigma_n = \mathcal{I}_n^{-1/2}\mathcal{I}_n(\psi_n,
\hat\lambda_n)\mathcal{I}_n^{-1/2}$, so that $\|\Sigma_n - I_d\| \leq
\delta_n^M < 1$ for all large $n$. Then $\mathcal{I}_n(\psi_n, \hat\lambda_n)
= \mathcal{I}_n^{1/2}\Sigma_n\mathcal{I}_n^{1/2}$ is positive definite,
because $\Sigma_n$ is symmetric, as $\mathcal{I}_n(\psi_n, \cdot)$ has
symmetric values, with eigenvalues in $(0, 2)$. Thus
$\bar{\mathcal{I}}_n(\psi_n, \hat\lambda_n)$ is defined and positive definite,
and the leading $d_1 \times d_1$ block of $\mathcal{I}_n(\psi_n,
\hat\lambda_n)^{-1}$ is $\bar{\mathcal{I}}_n(\psi_n, \hat\lambda_n)^{-1} =
[I_{d_1}, 0]\mathcal{I}_n^{-1/2}\Sigma_n^{-1}\mathcal{I}_n^{-1/2}[I_{d_1},
0]^\tsp$. Hence
\[
  \bar{\mathcal{I}}_n^{1/2}\bar{\mathcal{I}}_n(\psi_n, \hat\lambda_n)^{-1}\bar{\mathcal{I}}_n^{1/2} - I_{d_1}
  = L_n(\Sigma_n^{-1} - I_d)L_n^\tsp,
\]
whose norm is at most $\|\Sigma_n^{-1} - I_d\| \leq \delta_n^M/\{1 -
\delta_n^M\}$. This proves part (i).

For part (ii), let $Q_n = \mathcal{I}_{n\lambda\lambda}^{-1/2}[0,
I_{d_2}]\mathcal{I}_n^{1/2}$, which satisfies $Q_nQ_n^\tsp =
\mathcal{I}_{n\lambda\lambda}^{-1/2}\mathcal{I}_{n\lambda\lambda}
\mathcal{I}_{n\lambda\lambda}^{-1/2} = I_{d_2}$ and hence has spectral norm one, so
that $\tilde S_{n\lambda} = Q_n\tilde S_n(\theta_n)$ satisfies
$\|\tilde S_{n\lambda}\| \leq \|\tilde S_n(\theta_n)\|$. Fix $\eta \in (0, 1)$.
By Assumption \ref{assm:clt} there is an $M > 0$ with
$\pr\{2\|\tilde S_n(\theta_n)\| > M\} \leq \eta$ for all large $n$; fix such
an $M$ and let $\mathcal{G}_n$ be the event that $2\|\tilde S_{n\lambda}\| \leq
M$ and $\veps_n^M \leq 1/2$, so that $\pr(\mathcal{G}_n) \geq 1 - 2\eta$ for
all large $n$. Since $\eta$ is
arbitrary, it suffices to show that, on $\mathcal{G}_n$ and for all large
$n$, $S_{n\lambda}(\psi_n, \cdot)$ has exactly one root $\lambda^\ast$ in
$B_n^M$, that $\mathcal{J}_{n\lambda\lambda}(\psi_n, \lambda^\ast)$ is positive definite,
and that $\|\mathcal{I}_{n\lambda\lambda}^{1/2}(\lambda^\ast - \lambda_n)\| \leq
2\|\tilde S_{n\lambda}\|$. Indeed, in that case, let $\hat\lambda_n$ be the root of
$S_{n\lambda}(\psi_n, \cdot)$ with
$\|\mathcal{I}_{n\lambda\lambda}^{1/2}(\hat\lambda_n - \lambda_n)\| \leq
2\|\tilde S_{n\lambda}\|$ when exactly one such root exists, and
$\hat\lambda_n = \lambda_n$ otherwise, so that $\hat\lambda_n$ satisfies
\eqref{eq:info_bounded} because $2\|\tilde S_{n\lambda}\| = O_p(1)$. On
$\mathcal{G}_n$, a root with $\|\mathcal{I}_{n\lambda\lambda}^{1/2}(\lambda -
\lambda_n)\| \leq 2\|\tilde S_{n\lambda}\| \leq M$ is in $B_n^M$, so
$\hat\lambda_n = \lambda^\ast$, a root at which
$\mathcal{J}_{n\lambda\lambda}(\psi_n, \hat\lambda_n)$ is positive definite. If
$\hat\lambda_n'$ is a root satisfying \eqref{eq:info_bounded}, the $M$ above
can be taken so large that also $\pr(\hat\lambda_n' \notin B_n^M) \leq \eta$
for all large $n$, and on $\mathcal{G}_n \cap \{\hat\lambda_n' \in B_n^M\}$
uniqueness of the root in $B_n^M$ gives $\hat\lambda_n' = \hat\lambda_n$.

It remains to show that, on $\mathcal{G}_n$ and for all large $n$, there is
exactly one root $\lambda^\ast$ in $B_n^M$, with the two properties claimed.
For $u \in \R{d_2}$ with $\|u\| \leq M$, let $\lambda(u) =
\lambda_n + \mathcal{I}_{n\lambda\lambda}^{-1/2}u$, which is in $B_n^M$, and let
\[
  \Upsilon_n(u) = \tilde S_{n\lambda}\{\psi_n, \lambda(u)\}, \qquad
  \mathcal{K}_n(u) = \mathcal{I}_{n\lambda\lambda}^{-1/2}\mathcal{J}_{n\lambda\lambda}\{\psi_n,
  \lambda(u)\}\mathcal{I}_{n\lambda\lambda}^{-1/2},
\]
so that $\Upsilon_n(0) = \tilde S_{n\lambda}$, $\partial\Upsilon_n(u)/\partial u^\tsp =
-\mathcal{K}_n(u)$ by the chain rule and, for $\|u\|, \|u'\| \leq M$,
\begin{equation} \label{eq:rootdiff}
  \Upsilon_n(u) - \Upsilon_n(u') = -\int_0^1 \mathcal{K}_n\{u' + t(u - u')\}\,\dd t\,(u - u').
\end{equation}
Since $\mathcal{J}_{n\lambda\lambda}(\theta) = [0, I_{d_2}]\mathcal{J}_{n\cdot\lambda}(\theta)$,
$\mathcal{I}_{n\lambda\lambda} = [0, I_{d_2}]\mathcal{I}_{n\cdot\lambda}$ and
$Q_n\mathcal{I}_n^{-1/2} = \mathcal{I}_{n\lambda\lambda}^{-1/2}[0, I_{d_2}]$,
\[
  \mathcal{K}_n(u) - I_{d_2} = Q_n\mathcal{I}_n^{-1/2}\bigl[\mathcal{J}_{n\cdot\lambda}\{\psi_n, \lambda(u)\} -
  \mathcal{I}_{n\cdot\lambda}\bigr]\mathcal{I}_{n\lambda\lambda}^{-1/2},
\]
so $\sup_{\|u\| \leq M}\|\mathcal{K}_n(u) - I_{d_2}\| \leq \veps_n^M$
because $\|Q_n\| = 1$. Let $T_n(u) = u + \Upsilon_n(u)$, so that $u$ is a
root of $\Upsilon_n$ if and only if $T_n(u) = u$, and, by \eqref{eq:rootdiff},
$\|T_n(u) - T_n(u')\| \leq \veps_n^M\|u - u'\|$ for $\|u\|, \|u'\| \leq M$.
On $\mathcal{G}_n$, $T_n$ is a contraction with constant $1/2$ on the closed
ball of radius $M$ and maps it into itself, because $\|T_n(u)\| \leq
\|\tilde S_{n\lambda}\| + \veps_n^M\|u\| \leq M$ there. By Banach's
fixed-point theorem \citep[Theorem~5.1.2]{kreyszig1978introductory}, $T_n$
has exactly one fixed point $u_n^\ast$ in the ball, so $\lambda^\ast =
\lambda(u_n^\ast)$ is the only root of $S_{n\lambda}(\psi_n, \cdot)$ in $B_n^M$.
Moreover,
$\|u_n^\ast\| \leq \|\tilde S_{n\lambda}\| + \tfrac12\|u_n^\ast\|$, so
$\|\mathcal{I}_{n\lambda\lambda}^{1/2}(\lambda^\ast - \lambda_n)\| =
\|u_n^\ast\| \leq 2\|\tilde S_{n\lambda}\|$, and $\mathcal{J}_{n\lambda\lambda}(\psi_n,
\lambda^\ast) =
\mathcal{I}_{n\lambda\lambda}^{1/2}\mathcal{K}_n(u_n^\ast)\mathcal{I}_{n\lambda\lambda}^{1/2}$
is positive definite because $\|\mathcal{K}_n(u_n^\ast) - I_{d_2}\| \leq 1/2$. This completes the proof.
\end{proof}

\begin{proof}[Proof of Proposition \ref{prop:box}]
Fix an arbitrary sequence $(\theta_n)$ in $K$; it suffices to show that
every subsequence has a further subsequence along which $\pr\{\theta_n \in
C_n^{\mathrm{s}}(\alpha)\} \to 1 - \alpha$. Let $Z_n = \tilde S_n(\theta_n)$ and let
\[
  \mathcal{D}_n = \diag\{(\mathcal{I}_n^{-1})_{11}, \dots,
  (\mathcal{I}_n^{-1})_{dd}\}, \qquad G_n =
  \mathcal{D}_n^{-1/2}\mathcal{I}_n^{-1/2},
\]
so that $G_n Z_n =
\mathcal{D}_n^{-1/2}\mathcal{I}_n^{-1}S_n(\theta_n)$ and $G_nG_n^\tsp
= \Omega_n(\theta_n)$. With the $j$th coordinate as the parameter of
interest, the block form of the inverse gives $(\mathcal{I}_n^{-1}S_n)_j =
(\mathcal{I}_n^{-1})_{jj}U_{n,j}$ for the efficient score $U_{n,j}$ of
\eqref{eq:eff_score}, so $\{(\mathcal{I}_n^{-1})_{jj}\}^{1/2}U_{n,j}(\theta_n)
= (G_n Z_n)_j$, and the proofs of Lemma \ref{lem:chi2} and Theorem
\ref{thm:prof} give $\hat{R}_{n,j}(\theta_{n,j}) = (G_n Z_n)_j + o_p(1)$
under Assumptions \ref{assm:eff_info_stable}--\ref{assm:nuis_est} for that
coordinate. Hence $\max_j |\hat{R}_{n,j}(\theta_{n,j})| = \|G_n Z_n\|_\infty +
o_p(1)$, and $\theta_n \in C_n^{\mathrm{s}}(\alpha)$ if and only if $\max_j
|\hat{R}_{n,j}(\theta_{n,j})| \leq m_\alpha(\hat{\Omega}_n)$.

Because $\|G_n\|^2 = \|\Omega_n(\theta_n)\| \leq d$, we may pass to a
subsequence along which $G_n \to G$, and then
$\Omega_n(\theta_n) \to \Omega = GG^\tsp$, a correlation matrix.
Then $G_n Z_n \weak G Z \sim \rN_d(0, \Omega)$ by Assumption
\ref{assm:clt} and Slutsky's theorem, and $\hat{\Omega}_n \to \Omega$ in
probability. The function $(\Omega, c) \mapsto \pr(\|\Omega^{1/2}Z\|_\infty
\leq c)$ is continuous and strictly increasing in $c$, so $m_\alpha$ is
continuous and $m_\alpha(\hat{\Omega}_n) \to m_\alpha(\Omega)$ in
probability. Since $\|\Omega^{1/2}Z\|_\infty$ has a continuous distribution
function, Slutsky's theorem gives
\[
\begin{aligned}
  \pr\{\theta_n \in C_n^{\mathrm{s}}(\alpha)\}
  &= \pr\{\|G_n Z_n\|_\infty + o_p(1) \leq m_\alpha(\hat{\Omega}_n)\} \\
  &\to \pr\{\|\Omega^{1/2}Z\|_\infty \leq m_\alpha(\Omega)\} = 1 - \alpha
\end{aligned}
\]
along the subsequence. The sequence was arbitrary, so Lemma \ref{lem:unif}(i) with $g_n(\theta) =
\pr_\theta\{\theta \in C_n^{\mathrm{s}}(\alpha)\}$ and $c = 1 - \alpha$ gives
the claim. This completes the proof.
\end{proof}


\subsection{Proofs for Section \ref{sec:irreg}}
\label{sec:irreg_proofs}

\begin{lemma}\label{lem:Akappa}
If Assumption \ref{assm:information} holds, then $\|\tilde\beta_n(\theta_n)\|^2 = \kappa_n^2/(1 -
\kappa_n^2)$ and, when $d_1 = d_2 = 1$, $\tilde\beta_n(\theta_n) = \kappa_n/(1 - \kappa_n^2)^{1/2}$.
\end{lemma}

\begin{proof}
By \eqref{eq:An} and the definition of $\beta_n$,
\[
  \tilde\beta_n\tilde\beta_n^\tsp = \bar{\mathcal{I}}_n^{-1/2}\mathcal{I}_{n\psi\lambda}
  \mathcal{I}_{n\lambda\lambda}^{-1}\mathcal{I}_{n\lambda\psi}\bar{\mathcal{I}}_n^{-1/2}
  = \bar{\mathcal{I}}_n^{-1/2}\mathcal{I}_{n\psi\psi}\bar{\mathcal{I}}_n^{-1/2} - I_{d_1},
\]
where the second equality uses
$\mathcal{I}_{n\psi\lambda}\mathcal{I}_{n\lambda\lambda}^{-1}\mathcal{I}_{n\lambda\psi}
= \mathcal{I}_{n\psi\psi} - \bar{\mathcal{I}}_n$. The eigenvalues of
$\bar{\mathcal{I}}_n^{-1/2}\mathcal{I}_{n\psi\psi}\bar{\mathcal{I}}_n^{-1/2}$
are those of $\bar{\mathcal{I}}_n^{-1}\mathcal{I}_{n\psi\psi}$, the
reciprocals of the eigenvalues of
$\mathcal{I}_{n\psi\psi}^{-1}\bar{\mathcal{I}}_n = I_{d_1} -
\mathcal{I}_{n\psi\psi}^{-1}\mathcal{I}_{n\psi\lambda}\mathcal{I}_{n\lambda\lambda}^{-1}\mathcal{I}_{n\lambda\psi}$.
The eigenvalues of the subtracted matrix are the squared canonical
correlations between $S_{n\psi}$ and $S_{n\lambda}$
\citep[Section 12.2]{anderson2003introduction}, the largest of which is
$\kappa_n^2$, so the smallest eigenvalue of
$\mathcal{I}_{n\psi\psi}^{-1}\bar{\mathcal{I}}_n$ is $1 - \kappa_n^2$.
Hence $\|\tilde\beta_n\|^2 = 1/(1 -
\kappa_n^2) - 1 = \kappa_n^2/(1 - \kappa_n^2)$. When $d_1 = d_2 = 1$, $\tilde\beta_n$
has the sign of $\mathcal{I}_{n\psi\lambda}$ by \eqref{eq:An}, which is the
sign of $\kappa_n$. This completes the proof.
\end{proof}

\begin{proof}[Proof of Lemma~\ref{lem:primitive}]
Substituting (i) into (ii) and multiplying by $\tilde\beta_n$, bounded by
Assumption~\ref{assm:bddcorr}, gives the display in Assumption~\ref{assm:proj}.
This completes the proof.
\end{proof}

For $C \subseteq \R{d_1}$, let $\partial C$ be its boundary and, for
$\veps > 0$, let $(\partial C)^\veps$ be the open set of points at distance
less than $\veps$ from $\partial C$. If $\|x - y\| < \veps$ and exactly one of
$x$ and $y$ is in $C$, the segment between them meets $\partial C$ at a point
within distance $\veps$ of both, so $x, y \in (\partial C)^\veps$. Hence, for
Borel $C$ and random vectors $X$ and $Y$ in $\R{d_1}$ on the same probability
space, splitting the event that $1_C(X) \neq 1_C(Y)$ according to whether
$\|X - Y\| < \veps$ gives
\begin{equation}\label{eq:segprob}
  |\pr(X \in C) - \pr(Y \in C)| \leq \pr\{1_C(X) \neq 1_C(Y)\}
  \leq \pr\{X \in (\partial C)^\veps\} + \pr(\|X - Y\| \geq \veps).
\end{equation}
The following lemma gives the constant $c_{d_1}$ used in Section \ref{sec:irreg}.

\begin{lemma}\label{lem:gauss}
With $c_{d_1}$ as defined before Proposition \ref{prop:bound}, for every
convex $C \subseteq \R{d_1}$, $\veps > 0$ and $v \in \R{d_1}$,
\begin{enumerate}
  \item[(i)~~] $\pr\{Z_\psi \in (\partial C)^\veps\} \leq 2c_{d_1}\veps$, and
  \item[(ii)~] $|\pr(Z_\psi \in C - v) - \pr(Z_\psi \in C)| \leq c_{d_1}\|v\|$
  if, in addition, $C$ is Borel.
\end{enumerate}
\end{lemma}

\begin{proof}
(i) The closure of $C$ is convex with the same boundary, so we may take
$C$ closed, and if $C$ is empty or all of $\R{d_1}$, then $\partial C$ is
empty. Otherwise, let $\delta_C(x)$ be the distance from $x$ to $C$ if
$x \notin C$ and minus the distance from $x$ to $\R{d_1} \setminus C$ if
$x \in C$, and let $C^t = \{x : \delta_C(x) \leq t\}$ for $t \in \R{}$.
Because $\partial C$ is in the closure of $C$ and in that of its complement,
$|\delta_C(x)|$ is at most the distance from $x$ to $\partial C$, so
$(\partial C)^\veps \subseteq (C^\veps \setminus C) \cup (C \setminus C^{-\veps})$.
By \citet[Theorem 1.2]{raic2019multivariate}, $c_{d_1}$ bounds the Gaussian
perimeter of every convex subset of $\R{d_1}$, so each of the two sets has
probability at most $c_{d_1}\veps$ under $Z_\psi$ by Proposition 1.1 of the
same paper, which gives (i).

(ii) The difference $\pr(Z_\psi \in C - v) - \pr(Z_\psi \in C)$ is zero if
$C$ is empty or all of $\R{d_1}$. Otherwise, it equals
$\pr\{Z_\psi \in (C - v)\setminus C\} - \pr\{Z_\psi \in C \setminus (C - v)\}$.
A point of $(C - v)\setminus C$ is outside $C$ and within distance $\|v\|$ of
it, and a point of $C \setminus (C - v)$ is in $C$ and within distance $\|v\|$
of its complement. With $C^t$ defined as in (i) but without taking $C$
closed, the two sets are therefore contained in $C^t \setminus C$ and
$C \setminus C^{-t}$ for every $t > \|v\|$. By the results cited in (i), each
has probability at most $c_{d_1}t$, so the difference of the two
probabilities is at most $c_{d_1}t$ in absolute value. Letting
$t \downarrow \|v\|$ gives (ii). This completes the proof.
\end{proof}

\begin{proof}[Proof of Theorem~\ref{thm:approx}]
Fix an arbitrary sequence $(\theta_n)$ in $K$. We show $g_n(\theta_n)
\to 0$ for
\[
  g_n(\theta) = \sup_{C \in \mathcal{C}}\ \bigl| \pr_\theta\{\hat{R}_n(\psi) \in C\}
  - \pr\{R^\ast_n(\theta) \in C\} \bigr| ;
\]
it suffices to show that every subsequence has a further subsequence along
which this holds. Write $R_n =
\bar{\mathcal{I}}_n^{-1/2} S_{n\psi}(\hat\theta_n)$, with $\hat\theta_n$
as in Section \ref{sec:limits}, so that $R_n = \tilde U_n +
\tilde\beta_n\{\tilde S_{n\lambda} - \Proj_n(\tilde S_{n\lambda})\} + o_p(1)$
by \eqref{eq:Vnlin}.

Let $W_n = (\tilde U_n, \tilde S_{n\lambda})$ and let $O_n$ be the $d \times
d$ matrix whose first $d_1$ rows are those of $\Gamma_n$ in the proof of
Lemma \ref{lem:chi2} and whose last $d_2$ rows are those of $Q_n$ in the proof
of Theorem \ref{thm:roots}, so that $W_n = O_n \tilde S_n$. Both blocks have
identity covariance matrix, and
\[
  \cov(\tilde U_n, \tilde S_{n\lambda}) =
  \bar{\mathcal{I}}_n^{-1/2}(\mathcal{I}_{n\psi\lambda} -
  \beta_n\mathcal{I}_{n\lambda\lambda})\mathcal{I}_{n\lambda\lambda}^{-1/2} = 0,
\]
so $O_nO_n^\tsp = \cov(W_n) = I_d$ and $O_n$ is orthogonal. Passing to a
subsequence, $O_n \to O$ with $O$ orthogonal, so $W_n \weak OZ \sim \rN_d(0,
I_d)$ by Assumption~\ref{assm:clt} and the continuous mapping theorem; that
is, $W_n \weak Z$. Define $f_n : \R{d} \to \R{d_1}$ by
\[
  f_n(z) = z_\psi + \tilde\beta_n\{z_\lambda - \Proj_n(z_\lambda)\},
\]
where $z_\psi$ and $z_\lambda$ are the first $d_1$ and the last $d_2$
coordinates of $z$, so that $R_n = f_n(W_n) + o_p(1)$ by \eqref{eq:Vnlin}
and $f_n(Z)$ has the distribution of $R^\ast_n$ in \eqref{eq:Tstar}. The set
$\mathcal H_n$ is closed and convex by Assumption~\ref{assm:proj} and
contains $0$, so $\Proj_n$ is single-valued, $\|\Proj_n(z) - \Proj_n(z')\|
\leq \|z - z'\|$ \citep[Theorem~3.16 and Proposition~4.16]{bauschke2017convex}
and
$\|z_\lambda - \Proj_n(z_\lambda)\| \leq \|z_\lambda\|$. Hence $\|f_n(z) -
f_n(z')\| \leq L\|z - z'\|$ and $\|f_n(z)\| \leq L\|z\|$ for all $z$ and
$z'$, with $L = 1 + 2\sup_n\|\tilde\beta_n\|$, which is finite by
Assumption~\ref{assm:bddcorr}. Because $(W_n)$ converges in distribution it
is uniformly tight, and so is $\{f_n(W_n)\}$, so $R_n = O_p(1)$.
Assumption~\ref{assm:eff_info_stable} and Lemma \ref{lem:sqrt} give
$\bar{\mathcal{I}}_n(\hat\theta_n)^{-1/2}\bar{\mathcal{I}}_n^{1/2} - I_{d_1} =
o_p(1)$, so $\hat R_n(\psi_n) =
\bar{\mathcal{I}}_n(\hat\theta_n)^{-1/2}\bar{\mathcal{I}}_n^{1/2} R_n = R_n +
o_p(1)$. Thus $\hat R_n(\psi_n) = f_n(W_n) + \rho_n$ with $\rho_n = o_p(1)$.

Because $W_n \weak Z$ and $\rho_n = o_p(1)$, $(W_n, \rho_n) \weak (Z, 0)$
for any extension of $\rho_n$ to a random variable
\citep[Theorem~2.7(v)]{vandervaart1998asymptotic}, so by the almost-sure
representation theorem \citep[Theorem~2.19]{vandervaart1998asymptotic} there are
$(\tilde W_n, \tilde\rho_n)$ and $\tilde Z$ on one probability space,
distributed as $(W_n, \rho_n)$ and $Z$, with $(\tilde W_n, \tilde\rho_n) \to
(\tilde Z, 0)$ almost surely. Then $\tilde R_n = f_n(\tilde W_n) +
\tilde\rho_n$ has the distribution of $\hat R_n(\psi_n)$, $f_n(\tilde Z)$ has
that of $R^\ast_n$, and
\[
  \Delta_n = \|\tilde R_n - f_n(\tilde Z)\| \leq L\|\tilde W_n - \tilde Z\|
  + \|\tilde\rho_n\| \to 0
\]
almost surely, so $\pr(\Delta_n \geq \veps) \to 0$ for every $\veps > 0$, and
this probability does not depend on $C$.

Fix $\veps > 0$ and $C \in \mathcal{C}$. Applying \eqref{eq:segprob} with
$X = f_n(\tilde Z)$ and $Y = \tilde R_n$ gives
\[
\begin{aligned}
  \bigl|\pr\{\hat R_n(\psi_n) \in C\} - \pr(R^\ast_n \in C)\bigr|
  &= \bigl|\pr(\tilde R_n \in C) - \pr\{f_n(\tilde Z) \in C\}\bigr| \\
  &\leq \pr\{R^\ast_n \in (\partial C)^\veps\} + \pr(\Delta_n \geq \veps).
\end{aligned}
\]
Conditionally on $Z_\lambda$, $R^\ast_n = Z_\psi + w$ with $w =
\tilde\beta_n\{Z_\lambda - \Proj_n(Z_\lambda)\}$ fixed and $Z_\psi$ standard
normal, and $(\partial C)^\veps - w = \{\partial(C - w)\}^\veps$ with $C - w$
convex, so $\pr\{R^\ast_n \in (\partial C)^\veps \mid Z_\lambda\} \leq
2c_{d_1}\veps$ by Lemma \ref{lem:gauss}(i), and hence
$\pr\{R^\ast_n \in (\partial C)^\veps\} \leq 2c_{d_1}\veps$. Therefore
$\sup_{C \in \mathcal{C}} |\pr\{\hat R_n(\psi_n) \in C\} - \pr(R^\ast_n \in
C)| \leq 2c_{d_1}\veps + o(1)$. Letting $n \to \infty$ and then $\veps \to 0$
gives $g_n(\theta_n) \to 0$ along the subsequence.
The sequence was arbitrary, so Lemma \ref{lem:unif}(i) with $c = 0$ gives the
theorem. This completes the proof.
\end{proof}

\begin{proof}[Proof of Corollary \ref{cor:naive}]
For $\alpha \in (0, 1)$, let $B_\alpha$ be the closed ball in $\R{d_1}$ of
radius $q_{d_1, 1-\alpha}^{1/2}$ about the origin, which is convex. Then
$\{\psi \in C_n(\alpha)\} = \{\hat R_n(\psi) \in B_\alpha\}$ and
$\{\|R^\ast_n(\theta)\|^2 \leq q_{d_1, 1-\alpha}\} = \{R^\ast_n(\theta) \in
B_\alpha\}$, so the supremum in the claim is at most the supremum in Theorem
\ref{thm:approx}, which tends to zero. This completes the proof.
\end{proof}

\begin{proof}[Proof of Proposition \ref{prop:bound}]
Fix $C \in \mathcal{C}$. Given $Z_\lambda$, the vector $R^\ast_n$ is normal
with mean $w = \tilde\beta_n\{Z_\lambda - \Proj_n(Z_\lambda)\}$ and identity
covariance matrix, so Lemma \ref{lem:gauss}(ii) gives $|\pr(R^\ast_n \in C
\mid Z_\lambda) - \pr(Z_\psi \in C)| \leq c_{d_1}\|w\| \leq
c_{d_1}\|\tilde\beta_n\|\,\|Z_\lambda - \Proj_n(Z_\lambda)\|$. Since $r_n =
\inf\{\|z\| : z \notin \mathcal{H}_n\}$ and $\mathcal{H}_n$ is closed,
$\mathcal{H}_n$ contains the closed ball of radius $r_n$ about the origin,
so $\|Z_\lambda - \Proj_n(Z_\lambda)\| \leq \max(\|Z_\lambda\| - r_n, 0)$.
Taking expectations and then the supremum over $C$ gives the claim. This
completes the proof.
\end{proof}

\begin{proof}[Proof of Proposition \ref{prop:lf}]
For $|\kappa| < 1$, write $b(\kappa) = \kappa(1 - \kappa^2)^{-1/2}$, let
$R^\ast_{b,h}$ be as in Section \ref{sec:irreg} for $b \in \R{}$ and
$h \in [0, \infty)$, and let $R^\ast_{b,\infty} = Z_\psi$. By Lemma \ref{lem:Akappa},
$\tilde\beta_n(\theta_n) = b(\kappa_n)$, and $\Proj_n(Z_\lambda) =
\max(Z_\lambda, -h_n)$ when $Z_\lambda \leq u_n$, so $R^\ast_n =
R^\ast_{b(\kappa_n), h_n}$ on $\{Z_\lambda \leq u_n\}$.
For
$b, t \in \R{}$, let $G_b(t) = \min[\E\,\Phi\{t - b\min(Z_\lambda, 0)\},
\Phi(t)]$, the smaller of the distribution functions of $R^\ast_{b,0}$ and
$Z_\psi$ at $t$. Because $R^\ast_{b,0} \geq Z_\psi$ if $b < 0$ and
$R^\ast_{b,0} \leq Z_\psi$ if $b \geq 0$, $G_b$ is the distribution function
of $R^\ast_{b,0}$ in the first case and $\Phi$ in the second, so $G_b$ is
continuous and strictly increasing, and $c_\alpha(\kappa)$ is the $(1 -
\alpha)$-quantile of $G_{b(\kappa)}$. Because $\Phi$ is Lipschitz with
constant $\phi(0)$, where $\phi$ is the standard normal density, and
$\E|\min(Z_\lambda, 0)| = \phi(0)$, $|G_b(t) -
G_{b'}(t)| \leq \phi(0)^2 |b - b'|$ for all $b$, $b'$ and $t$. Finally,
$\pr(R^\ast_{b,h} > t) \leq 1 - G_b(t)$ for every $h \geq 0$ and $t$, because
$R^\ast_{b,h} \leq R^\ast_{b,0}$ if $b < 0$ and $R^\ast_{b,h} \leq Z_\psi$ if
$b \geq 0$.

Fix an arbitrary sequence $(\theta_n)$ in $K$, write $b_n = b(\kappa_n)$ and
$\hat b_n = b\{\kappa_n(\hat\theta_n)\}$, and let $\eta_n$ be the supremum in
Theorem \ref{thm:approx}, so that $\eta_n \to 0$. By Assumption
\ref{assm:bddcorr}, $\kappa_n$ is eventually in a compact subset of $(-1, 1)$,
on which $b$ is uniformly continuous, so $\hat b_n - b_n = o_p(1)$.
Fix $\veps > 0$ and $\alpha \in (0, 1)$. If $\alpha + \phi(0)^2\veps \geq
1$, the bound below holds because its left side is at most $1 - \alpha \leq
\phi(0)^2\veps$.
Now consider $\alpha + \phi(0)^2\veps < 1$. Since
$c_\alpha\{\kappa_n(\hat\theta_n)\}$ is the $(1 - \alpha)$-quantile of the
continuous and strictly increasing $G_{\hat b_n}$, the event
$[\hat{R}_n(\psi_n) > c_\alpha\{\kappa_n(\hat\theta_n)\}]$ is $\{G_{\hat
b_n}(\hat{R}_n(\psi_n)) > 1 - \alpha\}$, which on $\{|\hat b_n - b_n| \leq
\veps\}$ is contained in $\{G_{b_n}(\hat{R}_n(\psi_n)) > 1 - \alpha -
\phi(0)^2\veps\} = \{\hat{R}_n(\psi_n) > t_n\}$, with $t_n$ the $(1 - \alpha -
\phi(0)^2\veps)$-quantile of $G_{b_n}$. Hence,
\[
\begin{aligned}
  \pr[\hat{R}_n(\psi_n) > c_\alpha\{\kappa_n(\hat\theta_n)\}]
  &\leq \pr\{\hat{R}_n(\psi_n) > t_n\} + \pr(|\hat b_n - b_n| > \veps) \\
  &\leq \pr(R^\ast_n > t_n) + \pr(|\hat b_n - b_n| > \veps) + \eta_n,
\end{aligned}
\]
where the second inequality holds by the definition of $\eta_n$, because half-lines are convex. Since $R^\ast_n=R^\ast_{b_n,h_n}$ on $\{Z_\lambda\le u_n\}$, which has probability $\Phi(u_n)$, and $h_n\ge0$, $\pr(R^\ast_n > t_n) \leq 1 -
G_{b_n}(t_n) + \Phi(-u_n) = \alpha + \phi(0)^2\veps + \Phi(-u_n)$. Therefore
\[
\begin{aligned}
  &\sup_{\alpha \in (0, 1)} \bigl( \pr[\hat{R}_n(\psi_n) > c_\alpha\{\kappa_n(\hat\theta_n)\}]
  - \alpha \bigr) \\
  &\quad \leq \phi(0)^2\veps + \Phi(-u_n) + \pr(|\hat b_n - b_n| > \veps) + \eta_n,
\end{aligned}
\]
where the right-hand side tends to $\phi(0)^2\veps$ by the hypotheses.
Letting $\veps \to 0$ gives $\limsup_{n \to \infty} g_n(\theta_n) \leq 0$ for
$g_n(\theta) = \sup_{\alpha \in (0, 1)} (\pr_\theta[\hat{R}_n(\psi) >
c_\alpha\{\kappa_n(\hat\theta_n)\}] - \alpha)$. The sequence was
arbitrary, so Lemma \ref{lem:unif}(ii) with $c = 0$ gives the first claim.

For the second claim, fix $\alpha \in (0, 1)$ and a sequence $(\theta_n)$ in
$K$ along which $h_n \to 0$ and $\kappa_n \to \bar\kappa < 0$, and write $\bar
b = b(\bar\kappa)$. Since $G_{b'}(t)$ is continuous in $(b', t)$ and strictly
increasing in $t$, $c_\alpha$ is continuous, so $c_\alpha\{\kappa_n(\hat\theta_n)\}
- c_\alpha(\kappa_n) = o_p(1)$ and $c_\alpha(\kappa_n) \to
c_\alpha(\bar\kappa)$. Fix $\veps > 0$. Then
\[
\begin{aligned}
  \pr[\hat{R}_n(\psi_n) > c_\alpha\{\kappa_n(\hat\theta_n)\}]
  &\geq \pr\{\hat{R}_n(\psi_n) > c_\alpha(\kappa_n) + \veps\} - o(1) \\
  &\geq \pr\{R^\ast_n > c_\alpha(\kappa_n) + \veps\} - \eta_n - o(1),
\end{aligned}
\]
and, by the same comparison, $\pr\{R^\ast_n > c_\alpha(\kappa_n) + \veps\}
\geq 1 - \E\,\Phi\{c_\alpha(\kappa_n) + \veps - b_n\min(Z_\lambda + h_n, 0)\}
- \Phi(-u_n)$. By dominated convergence the right-hand side tends to $1 -
\E\,\Phi\{c_\alpha(\bar\kappa) + \veps - \bar b\min(Z_\lambda, 0)\} =
\pr\{R^\ast_{\bar b,0} > c_\alpha(\bar\kappa) + \veps\}$, which tends to
$\alpha$ as $\veps \to 0$ because $c_\alpha(\bar\kappa)$ is the $(1 -
\alpha)$-quantile of the continuous distribution of $R^\ast_{\bar b,0}$ when
$\bar\kappa < 0$. Hence the limit inferior of $\sup_{\theta \in K}
\pr_\theta[\hat{R}_n(\psi) > c_\alpha\{\kappa_n(\hat\theta_n)\}]$ is at least
$\alpha$, and with the first claim its limit superior is $\alpha$. This
completes the proof.
\end{proof}


\subsection{Proofs for Section \ref{sec:zip}}
\label{sec:zip_proofs}

Throughout this subsection, let
\[
  b(\psi) = e^{-\psi}, \qquad
  p_0(\theta) = \lambda + (1-\lambda)b(\psi) = b + \lambda(1-b) ,
\]
suppress the arguments of $b$ and $p_0$, and write $\hat p_{0,n} = n_0/n$ for
the sample proportion of zeros, where $n_0$ is the number of zeros, and let
$\mathcal{E}_n = \{0 < n_0 < n\}$, the event that the sample contains a zero and a
positive observation. Here $b$
is the probability that a $\mathrm{Poisson}(\psi)$ variable is zero, $p_0 =
\pr_\theta(Y = 0)$, and $1 - p_0 = (1-b)(1-\lambda)$. The endpoints of
$\tilde\Lambda(\psi)$ in \eqref{eq:zipext} are the values at which $p_0 = 0$
and $\lambda = 1$, so $p_0 \in (0,1)$ and $1 - \lambda > 0$ on
$\tilde{\Theta}$. The scores for one observation are
\begin{equation}
\label{eq:zipscore}
\begin{aligned}
  S_{\psi}(y;\theta) &= -\mathbf 1\{y=0\}\frac{(1-\lambda)b}{p_0}
  + \mathbf 1\{y>0\}\Bigl(\frac{y}{\psi}-1\Bigr), \\
  S_{\lambda}(y;\theta) &= \mathbf 1\{y=0\}\frac{1-b}{p_0}
  - \mathbf 1\{y>0\}\frac{1}{1-\lambda},
\end{aligned}
\end{equation}
and the information for one observation is
\begin{equation}
\label{eq:zipinfo}
\mathcal{I}(\theta)=
\begin{pmatrix}
\dfrac{(1-\lambda)^2b^2}{p_0} + \dfrac{1-\lambda}{\psi} + \lambda - p_0
& -\dfrac{b}{p_0}\\[8pt]
-\dfrac{b}{p_0}
& \dfrac{(1-b)^2}{p_0(1-p_0)}
\end{pmatrix}.
\end{equation}

The next lemma, used in the proofs of Lemma \ref{lem:zipinfo} and Theorem
\ref{thm:zip}, justifies differentiating under the expectation on the extended
set.

\begin{lemma}
\label{lem:zipreg}
If $\theta \in \tilde\Theta$, then $\E_\theta\{S(Y;\theta)\} = 0$ and
$\E_\theta\bigl\{ \partial S(Y;\theta)/\partial\theta^\tsp \bigr\} =
-\mathcal{I}(\theta)$.
\end{lemma}

\begin{proof}
By \eqref{eq:zipmodel}, the mass function $p(y;\theta) = \pr_\theta(Y = y)$
of Section \ref{sec:zip} equals $\lambda\mathbf 1\{y = 0\} + (1-\lambda)q(y;\psi)$,
where $q(y;\psi) = e^{-\psi}\psi^y/y!$ is the $\mathrm{Poisson}(\psi)$ mass
function, and it is strictly positive on $\tilde\Theta$ because there, $p_0 > 0$ and $1-\lambda > 0$, and thus $p(0; \theta) = p_0 > 0$ and
$p(k;\theta) = (1 - \lambda)q(k;\psi) > 0$ for $k \geq 1$. The derivatives $\partial q(y;\psi)/\partial\psi=q(y;\psi)(y/\psi-1)$ and $\partial^2q(y;\psi)/\partial\psi^2=q(y;\psi)\{(y/\psi-1)^2-y/\psi^2\}$ are $q(y;\psi)$ times a function of $y$, so their sums over $y$ are the expectations of $Y/\psi-1$ and of $(Y/\psi-1)^2-Y/\psi^2$ under $Y \sim \mathrm{Poisson}(\psi)$, with both sums being zero by $\E Y = \E (Y - \psi)^2 = \psi$.
Because
$p(y;\theta)$ is affine in $\lambda$ and $\sum_y q(y;\psi) = 1$, it follows
that
\[
  \sum_y \frac{\partial p(y;\theta)}{\partial\theta} = 0, \qquad
  \sum_y \frac{\partial^2 p(y;\theta)}{\partial\theta\partial\theta^\tsp} = 0.
\]
Since $S(y;\theta) = \{\partial p(y;\theta)/\partial\theta\}/p(y;\theta)$,
the first identity is $\E_\theta\{S(Y;\theta)\} = 0$, the first claim,
and hence $\E_\theta\{S(Y;\theta)S(Y;\theta)^\tsp\} =
\mathcal{I}(\theta)$. Taking expectations under $\theta$ in
\[
  \frac{\partial S(y;\theta)}{\partial\theta^\tsp}
  = \frac{1}{p(y;\theta)}\frac{\partial^2 p(y;\theta)}{\partial\theta\partial\theta^\tsp}
  - S(y;\theta)S(y;\theta)^\tsp
\]
and using the second identity gives the second claim. This completes the
proof.
\end{proof}

\begin{lemma}
\label{lem:zipinfo}
The map $\theta \mapsto \mathcal{I}(\theta)$ is continuous and positive
definite on $\tilde{\Theta}$, and $\theta \mapsto \bar{\mathcal{I}}(\theta)^{-1}$
is continuous and strictly positive on $\tilde{\Theta}$.
\end{lemma}

\begin{proof}
Every entry of \eqref{eq:zipinfo} is continuous on $\tilde\Theta$ because
$\psi$, $p_0$ and $1-p_0$ are strictly positive there. For positive
definiteness, fix $\theta \in \tilde{\Theta}$ and $u = (u_1,u_2)^\tsp \in
\R{2}$ with $u^\tsp \mathcal{I}(\theta) u = 0$. Because the score has mean
zero by Lemma \ref{lem:zipreg}, $0 = u^\tsp \mathcal{I}(\theta) u = \E_\theta[\{u^\tsp
S(Y;\theta)\}^2]$, so $u^\tsp S(Y;\theta) = 0$ almost surely. Every $k \geq
1$ has $\pr_\theta(Y = k) > 0$ because $\psi > 0$ and $1 - \lambda > 0$, so by
\eqref{eq:zipscore}, $u_1(k/\psi - 1) = u_2/(1-\lambda)$ for every $k \geq
1$. Subtracting the equation for $k=1$ from that for $k=2$ gives $u_1/\psi =
0$, hence $u_1 = 0$ and then $u_2 = 0$. Finally, $\bar{\mathcal{I}} =
\mathcal{I}_{\psi\psi} - \mathcal{I}_{\psi\lambda}^2/\mathcal{I}_{\lambda\lambda}$
is continuous and strictly positive wherever $\mathcal{I}$ is continuous and
positive definite, and hence so is $\bar{\mathcal{I}}^{-1}$. This completes
the proof.
\end{proof}

\begin{lemma}
\label{lem:zipnuis}
If $K \subseteq \Theta$ is compact, $(\theta_n)$ is a sequence in $K$, and
$\hat\lambda_n$ is as in \eqref{eq:zipnuisdef}, then
$\pr(\mathcal{E}_n) \to 1$, $S_{n\lambda}(\psi_n, \hat\lambda_n) = 0$ on
$\mathcal{E}_n$, and $\hat\lambda_n - \lambda_n = O_p(n^{-1/2})$.
\end{lemma}

\begin{proof}
Write $p_{0,n} = p_0(\theta_n)$. Because $p_0$ is continuous with values in
$(0,1)$ on $\Theta$, there is an $r < 1$, depending only on $K$, with $p_0
\leq r$ and $1 - p_0 \leq r$ on $K$. Hence
\[
  \pr(\mathcal{E}_n^c) = p_{0,n}^{\,n} + (1-p_{0,n})^n \leq 2r^n \to 0 .
\]

For the other two conclusions, fix $\psi > 0$. By \eqref{eq:zipmodel} and $1
- \lambda = (1-p_0)/(1-b)$,
\[
\ell_n(\psi, \lambda) = n \bigl\{ \hat{p}_{0,n} \log p_0 + (1-\hat{p}_{0,n}) \log (1-p_0) \bigr\} + C,
\]
where $C$ does not depend on $\lambda$. On $\mathcal{E}_n$, the right-hand
side is strictly concave in $p_0 \in (0, 1)$ and maximized uniquely at $p_0 =
\hat{p}_{0,n}$. The map $\lambda \mapsto p_0(\psi, \lambda)$ is an increasing
affine bijection from $\tilde\Lambda(\psi)$ onto $(0, 1)$, so $\lambda
\mapsto \ell_n(\psi, \lambda)$ is strictly concave on $\tilde\Lambda(\psi)$,
with unique maximizer
\begin{equation}
\label{eq:zippihat}
  \hat\lambda_n(\psi) = \frac{\hat p_{0,n} - b(\psi)}{1 - b(\psi)}
  \ \in\ \tilde\Lambda(\psi).
\end{equation}
Because $\hat{p}_{0,n} < 1$ on $\mathcal{E}_n$, \eqref{eq:zippihat} gives
$\hat\lambda_n(\psi) < 1$.

Summing the second component of \eqref{eq:zipscore} over the observations and
using $1/(1-\lambda) = (1-b)/(1-p_0)$,
\[
  S_{n\lambda}(\psi,\lambda)
  = n_0\frac{1-b}{p_0} - (n-n_0)\frac{1-b}{1-p_0}
  = n(1-b)\,\frac{\hat p_{0,n} - p_0}{p_0(1-p_0)} .
\]
Since $\hat p_{0,n} - p_0 = (1-b)\{\hat\lambda_n(\psi) - \lambda\}$ by
\eqref{eq:zippihat} and $\mathcal{I}_{n\lambda\lambda}(\psi,\lambda) =
n(1-b)^2/\{p_0(1-p_0)\}$ by \eqref{eq:zipsum} and \eqref{eq:zipinfo}, on
$\mathcal{E}_n$ and for every $\lambda \in \tilde\Lambda(\psi)$,
\begin{equation}
\label{eq:zipexact}
  S_{n\lambda}(\psi,\lambda)
  = \mathcal{I}_{n\lambda\lambda}(\psi,\lambda)\,
    \bigl\{\hat\lambda_n(\psi) - \lambda\bigr\} .
\end{equation}
With $\psi = \psi_n$ and $\lambda = \hat\lambda_n$, \eqref{eq:zipexact}
gives $S_{n\lambda}(\psi_n,\hat\lambda_n) = 0$ on $\mathcal{E}_n$. For the
rate, \eqref{eq:zippihat} gives $\hat\lambda_n - \lambda_n = (\hat p_{0,n}
- p_{0,n})/\{1-b(\psi_n)\}$, whose numerator has variance $p_{0,n}(1 -
p_{0,n})/n \leq 1/(4n)$ and whose denominator is bounded away from zero on
$K$, so Chebyshev's inequality gives $\hat\lambda_n - \lambda_n =
O_p(n^{-1/2})$. This completes the proof.
\end{proof}

\begin{lemma}
\label{lem:zipunif}
If $K \subseteq \Theta$ is compact and $(\theta_n)$ is a sequence in $K$,
then there are $\varrho > 0$ and $C < \infty$, depending only on $K$, such
that the following hold for every $n$.
\begin{enumerate}
\item[(i)~~] Every $\lambda$ with $|\lambda - \lambda_n| \leq \varrho$ is in
$\tilde\Lambda(\psi_n)$, and at every such $\lambda$ the first and second
derivatives in $\lambda$ of both components of $S(y;\psi_n,\lambda)$
are bounded in absolute value by $C$, for every $y \geq 0$.
\item[(ii)~] The score satisfies $\E\| \mathcal{I}(\theta_n)^{-1/2}
S(Y;\theta_n) \|^4 \leq C$.
\item[(iii)] The information satisfies $\|\mathcal{I}(\theta_n)^{-1/2}\| \leq
C$ and $C^{-1} \leq \mathcal{I}_{\lambda\lambda}(\theta_n) \leq C$.
\end{enumerate}
\end{lemma}

\begin{proof}
The functions $\psi$, $p_0$ and $1-\lambda$ are continuous and strictly
positive on the compact $K$, so each is bounded above and away from zero by
constants depending only on $K$. By Lemma \ref{lem:zipinfo}, $\mathcal{I}$
is continuous and positive definite on $\tilde\Theta \supseteq K$. The inverse
symmetric square root is continuous on positive definite matrices, so $\theta
\mapsto \mathcal{I}(\theta)^{-1/2}$ is continuous on $K$, and
$\mathcal{I}_{\lambda\lambda}$ is continuous and strictly positive there,
which gives (iii).

For (i), $p_0(\psi_n,\lambda) = p_0(\theta_n) + \{1-b(\psi_n)\}(\lambda -
\lambda_n)$, so taking $\varrho$ smaller than half the lower bounds for $p_0$
and $1 - \lambda$ on $K$ keeps $p_0(\psi_n,\lambda)$ and $1-\lambda$ above a
constant $\eta > 0$, depending only on $K$, whenever $|\lambda - \lambda_n|
\leq \varrho$, so such $\lambda$ are in $\tilde\Lambda(\psi_n)$.
Differentiating \eqref{eq:zipscore} in $\lambda$, using $\partial
p_0/\partial \lambda = 1-b$ and $p_0 + (1-\lambda)(1-b) = 1$,
\[
  \frac{\partial S_{\psi}}{\partial\lambda} = \mathbf 1\{y=0\}\frac{b}{p_0^2},
  \qquad
  \frac{\partial S_{\lambda}}{\partial\lambda}
  = -\mathbf 1\{y=0\}\frac{(1-b)^2}{p_0^2} - \mathbf 1\{y>0\}\frac{1}{(1-\lambda)^2},
\]
with second derivatives $-2\,\mathbf 1\{y=0\}b(1-b)/p_0^3$ and $2\,\mathbf
1\{y=0\}(1-b)^3/p_0^3 - 2\,\mathbf 1\{y>0\}(1-\lambda)^{-3}$. Each is bounded
in $y$ by a constant depending only on $\eta$, since $b \in (0,1)$.

For (ii), \eqref{eq:zipscore} gives $|S_{\lambda}(Y;\theta_n)| \leq
\max\bigl[\{1-b(\psi_n)\}/p_0(\theta_n),\ (1-\lambda_n)^{-1}\bigr]$ and
$|S_{\psi}(Y;\theta_n)| \leq b(\psi_n)/p_0(\theta_n) + Y/\psi_n + 1$, so by
the bounds of the first paragraph there are constants $C_1$ and $C_2$, free
of $n$, with $|S_{\lambda}(Y;\theta_n)| \leq C_1$ and $|S_{\psi}(Y;\theta_n)|
\leq C_2(1+Y)$.
Hence $\| S(Y ; \theta_n) \| \leq |S_\psi(Y;\theta_n)| + |S_\lambda(Y;\theta_n)| \leq (C_1 + C_2) (1+Y)$, and since $(1+Y)^4 \leq 8(1+Y^4)$,
\[
\| S(Y ; \theta_n) \|^4 \leq 8(C_1+C_2)^4(1+Y^4).
\]

Under $\theta_n$ the variable $Y$ is zero with probability
$\lambda_n$ and $\mathrm{Poisson}(\psi_n)$ otherwise, so $\E(Y^4)
= (1-\lambda_n)(\psi_n + 7\psi_n^2 + 6\psi_n^3 + \psi_n^4)$ is bounded on $K$,
and hence $\sup_n \E\|S(Y;\theta_n)\|^4 < \infty$, which with (iii) gives
(ii). This completes the proof.
\end{proof}

\begin{proof}[Proof of Theorem~\ref{thm:zip}]
Fix an arbitrary sequence $(\theta_n)$ in $K$. We verify the hypotheses of Theorem \ref{thm:roots}, with the extended set $\tilde\Lambda_n(\psi) =
\tilde\Lambda(\psi)$ of \eqref{eq:zipext}. Because $K$ is compact and
$\tilde\Theta$ is open, the set $K'$ of points within some distance
$\varrho' > 0$ of $K$ is a compact subset of $\tilde\Theta$. Let $\varrho$
and $C$ be as in Lemma \ref{lem:zipunif}.

By Lemma \ref{lem:zipinfo}, $\mathcal{I}_n(\theta_n) = n\,\mathcal{I}(\theta_n)$ is positive definite, giving Assumption \ref{assm:information}.

To show that Assumption \ref{assm:clt} holds, let $X_{n,i} = n^{-1/2}
\mathcal{I}(\theta_n)^{-1/2}S(Y_i;\theta_n)$, so that
$\tilde S_n(\theta_n) = \sum_{i=1}^n X_{n,i}$ by
\eqref{eq:zipsum}. These
vectors are independent with mean zero by Lemma \ref{lem:zipreg}, their
covariance matrices sum to
$I_2$, and $\sum_{i=1}^n \E\|X_{n,i}\|^4 \leq C/n$ by Lemma
\ref{lem:zipunif}(ii). Hence, for every $\veps > 0$,
\[
  \sum_{i=1}^n \E\bigl[\|X_{n,i}\|^2\mathbf 1\{\|X_{n,i}\| > \veps\}\bigr]
  \leq \veps^{-2}\sum_{i=1}^n \E\|X_{n,i}\|^4 \leq \frac{C}{\veps^2 n} \to 0,
\]
and \citet[Proposition 2.27]{vandervaart1998asymptotic} gives
$\tilde S_n(\theta_n) \weak \rN_2(0,I_2)$.

We now verify the remaining conditions of Theorem \ref{thm:roots}. By \eqref{eq:zipscore}, the scores are continuously differentiable in $\lambda$
on the open interval $\tilde\Lambda(\psi_n)$. For $\tilde r_n$, write $p(y;\theta) = \pr_\theta(Y = y)$, which by
\eqref{eq:zipmodel} is $\lambda\mathbf 1\{y = 0\} + (1-\lambda)q(y;\psi)$, where
$q(y;\psi) = e^{-\psi}\psi^y/y!$. Since $\partial p(y;\theta)/\partial\lambda = \mathbf 1\{y = 0\} - q(y;\psi)$, the nuisance score from one observation is $\{\mathbf 1\{y = 0\} - q(y;\psi)\}/p(y;\theta)$, and thus
\[
\mathcal{I}_{n\lambda\lambda}(\theta_n)
= \E\{S_{n\lambda}(\theta_n)^2\}
= n\sum_{y=0}^\infty \left\{ \frac{\mathbf 1\{y = 0\} - q(y;\psi_n)}{p(y;\theta_n)} \right\}^2 p(y;\theta_n).
\]
Writing $p(y;\theta)=q(y;\psi)+\lambda\{\mathbf 1\{y=0\}-q(y;\psi)\}$, it follows that
$p(y;\psi_n,\lambda) - p(y;\theta_n) = (\lambda - \lambda_n)\{\mathbf 1\{y = 0\} - q(y;\psi_n)\}$, so for every
$\lambda \in \R{}$,
\begin{equation}
\label{eq:lambdasum}
  (\lambda - \lambda_n)^2\mathcal{I}_{n\lambda\lambda}(\theta_n) = n\sum_{y=0}^\infty
  \frac{\{p(y;\psi_n, \lambda) - p(y;\theta_n)\}^2}{p(y;\theta_n)}.
\end{equation}
Let $\lambda \notin \tilde\Lambda(\psi_n)$. If $\lambda \le -e^{-\psi_n}/\{1-e^{-\psi_n}\}$, then $p(0;\psi_n,\lambda) = e^{-\psi_n} + \lambda\{1 - e^{-\psi_n}\} \le 0$, so $p(0;\psi_n,\lambda) - p(0;\theta_n) \le -p(0;\theta_n) < 0$ and the $y=0$ summand in \eqref{eq:lambdasum} is at least $p(0; \theta_n)$.
Hence,
$(\lambda-\lambda_n)^2\mathcal{I}_{n\lambda\lambda}(\theta_n) \ge n\,p(0;\theta_n)
= n\,\pr(Y = 0)$. Otherwise, if $\lambda \ge 1$, then $p(y;\psi_n,\lambda) = (1-\lambda)q(y;\psi_n) \le 0$ for
every $y \ge 1$, so
$p(y;\psi_n,\lambda) - p(y;\theta_n) \le -p(y;\theta_n) < 0$ and each
summand in \eqref{eq:lambdasum} with $y \ge 1$ is at least $p(y;\theta_n)$.
Thus, $(\lambda-\lambda_n)^2\mathcal{I}_{n\lambda\lambda}(\theta_n) \ge n\sum_{y \ge 1}
p(y;\theta_n) = n\,\pr(Y > 0)$. In both cases $(\lambda-\lambda_n)^2\mathcal{I}_{n\lambda\lambda}(\theta_n) \ge
n\min\{\pr(Y=0), \pr(Y>0)\}$. Taking the infimum over $\lambda \notin
\tilde\Lambda(\psi_n)$ gives $\tilde r_n^2 \ge n\min\{\pr(Y = 0),\ \pr(Y > 0)\}$,
which is at least a constant times $n$ on $K$ because $p_0$ is continuous with
values in $(0,1)$ on the compact $K$, so $\tilde r_n \to \infty$.

For \eqref{eq:info_unif}, fix $M > 0$. By Lemma \ref{lem:zipunif}(iii), every $\lambda \in B_n^M$ satisfies
$|\lambda - \lambda_n| \leq M(C/n)^{1/2} < \min(\varrho, \varrho')$ for large $n$. Consider such $n$.
Then every point of the segment between $\lambda_n$ and $\lambda$ is also within $\varrho$ of $\lambda_n$, so Lemma \ref{lem:zipunif}(i) applies on it, and $(\psi_n,
\lambda) \in K'$ for $\lambda \in B_n^M$, because $(\psi_n,\lambda)$ is within $|\lambda-\lambda_n| < \varrho'$ of $\theta_n \in K$.

Because $\mathcal{I}_n =
n\mathcal{I}$, we have
$$
\mathcal{I}_n(\theta_n)^{-1/2}\{\mathcal{I}_n(\psi_n,\lambda) - \mathcal{I}_n(\theta_n)\}\mathcal{I}_n(\theta_n)^{-1/2} = \mathcal{I}(\theta_n)^{-1/2}\{\mathcal{I}(\psi_n,\lambda) - \mathcal{I}(\theta_n)\}\mathcal{I}(\theta_n)^{-1/2},
$$
and since
$\|\mathcal{I}(\theta_n)^{-1/2}\| \leq C$ by Lemma \ref{lem:zipunif}(iii), the second
supremum in \eqref{eq:info_unif} is at most $C^2\sup_{\lambda \in
B_n^M}\|\mathcal{I}(\psi_n, \lambda) - \mathcal{I}(\theta_n)\|$.
The information $\mathcal{I}$ is continuous on $\tilde\Theta$ by Lemma \ref{lem:zipinfo} and hence uniformly continuous on the compact set $K' \subseteq \tilde\Theta$. Let $\omega(t) = \sup\{\|\mathcal{I}(\theta) - \mathcal{I}(\theta')\| : \theta,\theta' \in K', \|\theta - \theta'\| \le t\}$, which tends to zero as $t \downarrow 0$ by uniform continuity and depends only on $K'$. Since $\lambda \in B_n^M$, $(\psi_n,\lambda)$ and $\theta_n$ both lie in $K'$ and $\|(\psi_n,\lambda) - \theta_n\| = |\lambda - \lambda_n| \le M(C/n)^{1/2}$, we have $C^2\sup_{\lambda \in
B_n^M}\|\mathcal{I}(\psi_n, \lambda) - \mathcal{I}(\theta_n)\| \leq C^2 \omega\bigl(M(C/n)^{1/2}\bigr) \to 0$, giving the second condition in \eqref{eq:info_unif}.

For $y \ge 0$, let $\mathcal{J}_{\cdot\lambda}(y;\psi,\lambda) = -\partial S(y;\psi,\lambda)/\partial\lambda \in \R{2}$, whose entries are $-\partial S_\psi/\partial\lambda$ and $-\partial S_\lambda/\partial\lambda$, so that $\mathcal{J}_{n\cdot\lambda}(\psi,\lambda) = \sum_{i=1}^n \mathcal{J}_{\cdot\lambda}(Y_i;\psi,\lambda)$.
For the first condition in \eqref{eq:info_unif}, write
\[
  \mathcal{J}_{n\cdot\lambda}(\psi_n, \lambda) - \mathcal{I}_{n\cdot\lambda}(\theta_n)
  = \{\mathcal{J}_{n\cdot\lambda}(\psi_n, \lambda) - \mathcal{J}_{n\cdot\lambda}(\theta_n)\}
  + \{\mathcal{J}_{n\cdot\lambda}(\theta_n) - \mathcal{I}_{n\cdot\lambda}(\theta_n)\}.
\]
Each entry of $\mathcal{J}_{\cdot\lambda}(y;\psi_n,\cdot)$ is continuously differentiable on $\tilde\Lambda(\psi_n)$, and its derivative is at most $C$ in absolute value at every $\lambda$ with $|\lambda - \lambda_n| \le \varrho$ and for every $y \ge 0$ by Lemma \ref{lem:zipunif}(i).
For $\lambda \in B_n^M$, $|\lambda - \lambda_n| \le M(C/n)^{1/2} \le \varrho$, so the segment between $\lambda_n$ and $\lambda$ satisfies $|\lambda - \lambda_n| \leq \varrho$, and the mean value theorem gives that the absolute value of each entry of $\mathcal{J}_{\cdot\lambda}(y;\psi_n,\lambda) - \mathcal{J}_{\cdot\lambda}(y;\psi_n,\lambda_n)$ is at most $C|\lambda - \lambda_n|$, for every $y$. Summing over the $n$ observations, each entry of $\mathcal{J}_{n\cdot\lambda}(\psi_n,\lambda) - \mathcal{J}_{n\cdot\lambda}(\theta_n)$ is at most $Cn|\lambda-\lambda_n|$ in absolute value, and the norm of this vector is at most the sum of the absolute values of its two entries, so $\|\mathcal{J}_{n\cdot\lambda}(\psi_n,\lambda) - \mathcal{J}_{n\cdot\lambda}(\theta_n)\| \leq 2Cn|\lambda-\lambda_n|$, so that $|\lambda - \lambda_n| \le M(C/n)^{1/2}$ gives $2Cn|\lambda-\lambda_n| \le 2C^{3/2}Mn^{1/2}$, so the first quantity in brackets has norm at most $2C^{3/2}Mn^{1/2}$, which does not depend on $\lambda \in B_n^M$.

The second term does not depend on $\lambda$ and is a
sum of $n$ independent terms, each with mean zero by Lemma \ref{lem:zipreg}
and bounded by Lemma \ref{lem:zipunif}(i).
Write $\mathcal{J}_{n\cdot\lambda}(\theta_n) - \mathcal{I}_{n\cdot\lambda}(\theta_n) = \sum_{i=1}^n \xi_i$ with $\xi_i = \mathcal{J}_{\cdot\lambda}(Y_i;\theta_n) - \mathcal{I}_{\cdot\lambda}(\theta_n)$, where $\mathcal{I}_{\cdot\lambda}$ is the $\lambda$ column of $\mathcal{I}$. The $\xi_i$ are independent and identically distributed. By the second claim of Lemma \ref{lem:zipreg}, $\E\{-\partial S(Y;\theta_n)/\partial\theta^T\} = \mathcal{I}(\theta_n)$, so $\E\,\mathcal{J}_{\cdot\lambda}(Y;\theta_n) = \mathcal{I}_{\cdot\lambda}(\theta_n)$ and $\E\,\xi_i = 0$. By Lemma \ref{lem:zipunif}(i) at $\lambda = \lambda_n$, each entry of $\mathcal{J}_{\cdot\lambda}(y;\theta_n)$ is at most $C$ in absolute value for every $y$, hence so is each entry of its mean, so each entry of $\xi_i$ is at most $2C$ in absolute value. Each entry of $\sum_i \xi_i$ therefore has mean zero and variance at most $4C^2n$, and Chebyshev's inequality gives, for every $t > 0$ and every $n$,
\[
  \pr\left(\left|\sum_{i=1}^n \xi_{i}^{(r)}\right| > t\,n^{1/2}\right) \le \frac{4C^2}{t^2}, \qquad r = 1,2,
\]
where $\xi_i^{(r)}$ is the $r$th entry. The bound does not depend on $n$, so both entries, and hence the norm of $\sum_i\xi_i$, are $O_p(n^{1/2})$, so $n^{-1}\sum_i\xi_i = O_p(n^{-1/2}) = o_p(1)$. Since the second term does not depend on $\lambda$, the supremum over $B_n^M$ of the norm of the sum of the two terms is at most $2C^{3/2}Mn^{1/2} + \|\sum_i\xi_i\|$, so together with the first term, $\sup_{\lambda \in B_n^M}\|\mathcal{J}_{n\cdot\lambda}(\psi_n,\lambda) - \mathcal{I}_{n\cdot\lambda}(\theta_n)\| = O_p(n^{1/2})$.
By Lemma \ref{lem:zipunif}(iii), $\|\mathcal{I}_n(\theta_n)^{-1/2}\| = n^{-1/2}\|\mathcal{I}(\theta_n)^{-1/2}\| \le Cn^{-1/2}$ and, because $\mathcal{I}_{\lambda\lambda}(\theta_n) \ge C^{-1}$, $\mathcal{I}_{n\lambda\lambda}(\theta_n)^{-1/2} = n^{-1/2}\mathcal{I}_{\lambda\lambda}(\theta_n)^{-1/2} \le C^{1/2}n^{-1/2}$, so $\|\mathcal{I}_n(\theta_n)^{-1/2}\| \cdot \mathcal{I}_{n\lambda\lambda}(\theta_n)^{-1/2} \leq C^{3/2}n^{-1}$. Thus the first supremum in \eqref{eq:info_unif}, is at most $C^{3/2}n^{-1}\cdot O_p(n^{1/2}) = O_p(n^{-1/2})$, which is $o_p(1)$.

By Lemma \ref{lem:zipnuis}, $\hat\lambda_n$ is a root of
the nuisance score on $\mathcal{E}_n$, whose probability tends to one, so
Assumption \ref{assm:nuis_est} holds. Lemma \ref{lem:zipnuis} also gives
$\hat\lambda_n - \lambda_n = O_p(n^{-1/2})$, which is \eqref{eq:info_bounded}
because $\mathcal{I}_{n\lambda\lambda}(\theta_n) \leq Cn$ by Lemma
\ref{lem:zipunif}(iii) (so $\|\mathcal{I}_{n\lambda\lambda}(\theta_n)^{1/2}(\hat\lambda_n - \lambda_n)\| \le (Cn)^{1/2}\,O_p(n^{-1/2}) = O_p(1)$), so Theorem \ref{thm:roots}(i) gives Assumptions
\ref{assm:eff_info_stable} and \ref{assm:linear_approx}. The sequence was
arbitrary, so Theorem \ref{thm:prof} gives \eqref{eq:coverage}, and its proof
gives $\hat R_n(\psi_n) \rightsquigarrow \rN(0, 1)$. This completes the proof.
\end{proof}


\begin{proof}[Proof of Lemma \ref{lem:zipproj}]
Fix an arbitrary sequence $(\theta_n)$ in $K$. Assumptions \ref{assm:information} and \ref{assm:clt} hold as
in the proof of Theorem \ref{thm:zip}, where the estimator does not enter.
For Assumption \ref{assm:bddcorr}, $\tilde\beta_n(\theta)$ is continuous on
$\tilde\Theta$ by Lemma \ref{lem:zipinfo}, hence bounded on the compact $K$.
The set $\Lambda(\psi_n) = [0,\bar\lambda]$ is closed and convex,
$\mathcal{H}_n$ is the interval \eqref{eq:zipH}, and $u_n$ is at least a
positive constant times $n^{1/2}$ because $\mathcal{I}_{\lambda\lambda}$ and
$\bar\lambda - \lambda$ are continuous and strictly positive on $K$.

Let $\check\lambda_n$ be the maximizer \eqref{eq:zippihat} over the extended
set at $\psi = \psi_n$, so that \eqref{eq:zipexact} holds with
$\check\lambda_n$ in place of $\hat\lambda_n(\psi)$. On the event
$\mathcal{E}_n$ of Lemma \ref{lem:zipnuis}, $\lambda \mapsto
\ell_n(\psi_n,\lambda)$ is strictly concave on $\tilde\Lambda(\psi_n)
\supseteq \Lambda(\psi_n)$ with maximizer $\check\lambda_n$, and
\eqref{eq:zipexact} at $\lambda = \lambda_n$ gives $\check\lambda_n =
\lambda_n + \mathcal{I}_{n\lambda\lambda}(\theta_n)^{-1/2}\tilde S_{n\lambda}$. A strictly
concave function on an interval is maximized over a closed subinterval at the
point of the subinterval nearest its maximizer, so $\hat\lambda_n$ is the
point of $\Lambda(\psi_n)$ nearest $\check\lambda_n$, and
\[
 \mathcal{I}_{n\lambda\lambda}(\theta_n)^{1/2}\bigl(\hat\lambda_n - \lambda_n\bigr)
  = \Proj_n(\tilde S_{n\lambda}) ,
\]
which is condition (i) of Lemma \ref{lem:primitive}. Because $0 \in \mathcal
H_n$, $|\Proj_n(\tilde S_{n\lambda})| \leq |\tilde S_{n\lambda}| = O_p(1)$, so $\hat\lambda_n
- \lambda_n = O_p(n^{-1/2})$.

For condition (ii), apply \eqref{eq:zipexact}, with $\check\lambda_n$, at
$\lambda = \hat\lambda_n$ and at $\lambda = \lambda_n$, and divide by
$\mathcal{I}_{n\lambda\lambda}(\theta_n)^{1/2}$:
\[
\begin{aligned}
  \mathcal{I}_{n\lambda\lambda}(\theta_n)^{-1/2}
  \bigl\{ S_{n\lambda}(\psi_n,\hat\lambda_n) &- S_{n\lambda}(\theta_n) \bigr\} \\
  &= -\Proj_n(\tilde S_{n\lambda})
    + (s_n - 1)\bigl\{ \tilde S_{n\lambda} - \Proj_n(\tilde S_{n\lambda}) \bigr\},
\end{aligned}
\]
where $s_n =
\mathcal{I}_{n\lambda\lambda}(\psi_n,\hat\lambda_n)/\mathcal{I}_{n\lambda\lambda}(\theta_n)$.
Since $\hat\lambda_n - \lambda_n = O_p(n^{-1/2})$ and
$\mathcal{I}_{\lambda\lambda}$ is uniformly continuous and bounded away from
zero on the set $K'$ of the proof of Theorem \ref{thm:zip}, $s_n - 1 =
o_p(1)$, and $|\tilde S_{n\lambda} - \Proj_n(\tilde S_{n\lambda})| \leq |\tilde S_{n\lambda}| = O_p(1)$, so the
second term is $o_p(1)$ and Lemma \ref{lem:primitive} gives Assumption
\ref{assm:proj}. Finally, the rate $\hat\lambda_n - \lambda_n =
O_p(n^{-1/2})$ and Lemma \ref{lem:zipunif}(iii) give \eqref{eq:info_bounded},
so Theorem \ref{thm:roots}(i), whose hypotheses are verified in the proof of
Theorem \ref{thm:zip}, gives Assumptions \ref{assm:eff_info_stable} and
\ref{assm:linear_approx}. For Assumption \ref{assm:coef_stable}, its
left-hand side is $[\bar{\mathcal{I}}_n(\theta_n)^{-1/2}\{\beta_n(\hat\theta_n) -
\beta_n(\theta_n)\}\mathcal{I}_{n\lambda\lambda}(\theta_n)^{1/2}]\tilde
S_{n\lambda}(\hat\theta_n)$, where the first factor is $o_p(1)$ because
$\mathcal{I}$ is continuous on $\tilde\Theta$ by Lemma \ref{lem:zipinfo} and
$\hat\lambda_n - \lambda_n = O_p(n^{-1/2})$, and $\tilde S_{n\lambda}(\hat\theta_n)
= \tilde S_{n\lambda} - \Proj_n(\tilde S_{n\lambda}) + o_p(1) = O_p(1)$ by the
display above. This completes the proof.
\end{proof}


\subsection{Proofs for Section \ref{sec:oneway}}
\label{sec:oneway_proofs}

\begin{proof}[Proof of Proposition \ref{prop:oneway}]
Fix an arbitrary sequence $(\theta_n)$ in $K$; it suffices to show that
every subsequence has a further subsequence along which the claims hold.
Write $\tau_n = \tau(\theta_n)$ and let
\[
  \xi_n = \bigl[\{Q_{\mathrm{b}}/\tau_n - (m-1)\}/\{2(m-1)\}^{1/2},\ \{Q_{\mathrm{w}}/\psi_n - m(k-1)\}/\{2m(k-1)\}^{1/2}\bigr]^\tsp.
\]
Under $\theta_n$, $Q_{\mathrm{b}}/\tau_n$ and $Q_{\mathrm{w}}/\psi_n$ are independent $\chi^2_{m-1}$ and
$\chi^2_{m(k-1)}$ variables, so each component of $\xi_n$ is a standardized
sum of independent $\chi^2_1$ variables and $\xi_n \weak \rN_2(0, I_2)$ by the
central limit theorem and independence. By \eqref{eq:owscore},
$S_n(\theta_n) = \Xi_n\xi_n$ with
\[
  \Xi_n = \begin{bmatrix} \{(m-1)/2\}^{1/2}/\tau_n & \{m(k-1)/2\}^{1/2}/\psi_n \\
  k\{(m-1)/2\}^{1/2}/\tau_n & 0 \end{bmatrix},
\]
and $\mathcal{I}_n(\theta_n) = \Xi_n \Xi_n^\tsp$ by \eqref{eq:owinfo}. Since
$\Xi_n$ is invertible for $m \geq 2$, so is $\mathcal{I}_n(\theta_n)$, and
$\mathcal{I}_n(\theta_n)^{-1/2} \Xi_n$ is orthogonal.

Since $\mathcal{I}_n(\theta_n) = \Xi_n \Xi_n^\tsp$, the matrix
$\mathcal{I}_n(\theta_n)^{-1/2}\Xi_n$ satisfies
$\mathcal{I}_n^{-1/2}\,\Xi_n \Xi_n^\tsp\, \mathcal{I}_n^{-1/2} = I_2$,
so ir is orthogonal and all its entries lie in $[-1,1]$. The set of
$2\times 2$ orthogonal matrices is compact, so every subsequence has a
further subsequence along which $\mathcal{I}_n(\theta_n)^{-1/2}\Xi_n \to O$ with $O$
orthogonal.
Thus, $\tilde S_n(\theta_n) = \mathcal{I}_n(\theta_n)^{-1/2}\Xi_n\xi_n \weak OZ
\sim \rN_2(0, I_2)$, where $Z \sim \rN_2(0, I_2)$, by the continuous mapping
theorem, which is Assumption \ref{assm:clt}. Assumption \ref{assm:bddcorr}
holds because $\tilde\beta_n(\theta)^2 = (m-1)\psi^2/\{m(k-1) \tau^2\} \leq 1/(k-1)$ on
$\Theta$.

For Assumption \ref{assm:proj}, we verify conditions (i) and (ii) of Lemma
\ref{lem:primitive}. By \eqref{eq:owscore} and \eqref{eq:owinfo}, for every
$\lambda \in \tilde\Lambda(\psi_n)$,
\begin{equation} \label{eq:owexact}
  S_{n\lambda}(\psi_n, \lambda) = \mathcal{I}_{n\lambda\lambda}(\psi_n, \lambda)\bigl[\{Q_{\mathrm{b}}/(m-1)
  - \psi_n\}/k - \lambda\bigr],
\end{equation}
and $\tilde S_{n\lambda}$ is the first component of $\xi_n$. By
\eqref{eq:owexact} at $\lambda = \lambda_n$,
$\mathcal{I}_{n\lambda\lambda}(\theta_n)^{1/2} [\{Q_{\mathrm{b}}/(m-1) - \psi_n\}/k -
\lambda_n] = \tilde S_{n\lambda}$, so \eqref{eq:owconstr} gives
\[
  \mathcal{I}_{n\lambda\lambda}(\theta_n)^{1/2}(\hat\lambda_n - \lambda_n) = \max\{\tilde S_{n\lambda},
  -h_n(\theta_n)\} = \Proj_n(\tilde S_{n\lambda}),
\]
which is condition (i), because $\mathcal{H}_n = [-h_n(\theta_n), \infty)$.
Moreover, $\hat\tau_n = \psi_n + k\hat\lambda_n$ satisfies
$\hat\tau_n/\tau_n - 1 = \{2/(m-1)\}^{1/2}\Proj_n(\tilde S_{n\lambda})$ by
\eqref{eq:owinfo}, and $|\Proj_n(\tilde S_{n\lambda})| \leq |\tilde S_{n\lambda}|$ because $0
\in \mathcal{H}_n$, so $\hat\tau_n/\tau_n \to 1$ in probability. For
condition (ii), apply \eqref{eq:owexact} at $\lambda = \hat\lambda_n$ and at
$\lambda = \lambda_n$ and divide by
$\mathcal{I}_{n\lambda\lambda}(\theta_n)^{1/2}$:
\[
\begin{aligned}
  \mathcal{I}_{n\lambda\lambda}(\theta_n)^{-1/2}\{S_{n\lambda}(\psi_n, \hat\lambda_n) -
  S_{n\lambda}(\theta_n)\}
  &= (\tau_n/\hat\tau_n)^2\{\tilde S_{n\lambda} - \Proj_n(\tilde S_{n\lambda})\} - \tilde S_{n\lambda} \\
  &= -\Proj_n(\tilde S_{n\lambda}) + o_p(1),
\end{aligned}
\]
where the first equality uses $\mathcal{I}_{n\lambda\lambda}(\psi_n,
\hat\lambda_n)/\mathcal{I}_{n\lambda\lambda}(\theta_n) =
(\tau_n/\hat\tau_n)^2$ and the second that $\tilde S_{n\lambda} -
\Proj_n(\tilde S_{n\lambda}) = O_p(1)$. Lemma \ref{lem:primitive} gives
Assumption \ref{assm:proj}.

Assumptions \ref{assm:eff_info_stable} and \ref{assm:linear_approx} hold with
both left-hand sides equal to zero, because neither $U_n(\psi_n, \lambda)$
nor $\bar{\mathcal{I}}_n(\psi_n, \lambda)$ depends on $\lambda$,
and the quantity in Assumption \ref{assm:eff_info_stable} exists because
$\bar{\mathcal{I}}_n(\hat\theta_n) = \bar{\mathcal{I}}_n(\theta_n)$ and
$\mathcal{I}_{n\lambda\lambda}(\hat\theta_n) = k^2(m-1)/(2\hat\tau_n^2) > 0$
for $m \geq 2$ by \eqref{eq:owinfo}, as $\hat\tau_n \geq \psi_n > 0$.
Finally,
$\tilde\beta_n(\hat\theta_n) - \tilde\beta_n = (\tau_n/\hat\tau_n -
1)\tilde\beta_n = o_p(1)$. This completes the proof.
\end{proof}

\end{document}